\documentclass[12pt]{article} 

\usepackage{amssymb}
\usepackage{amsmath}
\usepackage{graphicx}
\usepackage{epsfig}
\usepackage{subfigure}

\def\shadowbox{\hbox{\rule[-0.0ex]{0.1ex}{1.2ex}%
\hspace{-0.1ex}\rule[-0.0ex]{1.2ex}{0.1ex}%
\hspace{0.0ex}\rule[-0.0ex]{0.1ex}{1.2ex}\hspace{-1.3ex}%
\rule[1.15ex]{1.25ex}{0.1ex}\hspace{-0.0ex}\rule[-0.25ex]{0.3ex}{1.1ex}%
\hspace{-1.2ex}\rule[-0.25ex]{1.1ex}{0.25ex}}}
\def\qed{\ifmmode \hbox{\hfill\shadowbox}
     \else \hphantom{x}\hfill\shadowbox \fi}

\newtheorem{theorem}{Theorem}[section]
\newtheorem{lemma}[theorem]{Lemma}
\newtheorem{definition}[theorem]{Definition}
\newtheorem{proposition}[theorem]{Proposition}
\newtheorem{corollary}[theorem]{Corollary}

\def\Cst{{\mathbb C}}

\def\Qst{{\mathbb Q}}
\def\Rst{{\mathbb R}}
\def\Tst{{\mathbb T}}
\def\Zst{{\mathbb Z}}

\def\zdd{\Zst ^{2d}}

\def\theta{\vartheta}
\def\phi{\varphi}
\def\Lsp{L}

\def\group{{\cal G}}
\def\hatg{\widehat{\group}}
\def\Ksub{{\cal K}}
\def\Kperp{{\cal K}^\perp}

\def\chiK{\chi_{\cal K}}
\def\lsp{\ell}
\def\Ltsp{{\Lsp^2}}

\def\ltZ{{\lsp^2(\Zst)}}
\def\ltsp{{\lsp^2}}
\def\llsp{{\lsp^2}(\Lambda)}

\def\rd{{\Rst^d}}

\def\toinf{{\rightarrow \infty }}

\def\Op{{\operatorname{Op}}}

\def\ran{{\operatorname{ran}}}
\def\ker{{\operatorname{ker}}}

\def\kn{K_{\sigma}}
\def\knt{K_{\tau}}
\def\knst{K_{\fouri(\hat{\sigma}  \twist   \hat{\tau})}}
\def\STFT{{\cal V}}

\def\tf{\group \times \hatg}
\def\cG{\mathcal{G}}

\def\cJ{\mathcal J}
\def\cJi{{\mathcal J}^{-1}}
\def\cS{\mathcal S}
\def\cC{\mathcal C}
\def\cA{\mathcal A}
\def\cF{\mathcal F}
\def\cV{\mathcal V}
\def\f{\phi}

\def\phas{(x,\omega)}
\def\vgf{V_g f}
\def\hattf{\hatg \times \group}

\def\mod{M^{\infty,1}}

\def\modv{M^{\infty,1}_v}
\def\modvj{M^{\infty,1}_{v\circ \cJ^{-1}}}
\def\modvij{M^{\infty,1}_{(v\otimes 1) \circ \cJ^{-1}}}
\def\modivj{M^{\infty,1}_{(1\otimes v) \circ \cJ^{-1}}}
\def\modii{M^{1}}
\def\modiim{M^{1}_m}
\def\modiiv{M^{1}_v}
\def\modiivi{M^{1}_{v\otimes 1}}
\def\modiiiv{M^{1}_{1\otimes v}}

\def\modiivg{M^{1}_v(\group)}

\def\specg{{\cal S}_{\cal C}(\group)}
\def\four{{\cal F}}
\def\fouri{{\cal F}^{-1}}
\def\twist{{\, \natural \, }}

\def\gabframe{\{{M_{\mu} T_{m} g}\}_{(m,\mu)\in\Lambda}}
\def\piframe{\{{\pi(\bfm ) g}\}_{\bfm \in\Lambda}}
\def\bfx{{\bf x}}
\def\bfy{{\bf y}}
\def\bfz{{\bf z}}
\def\bfu{{\bf u}}
\def\bfv{{\bf v}}

\def\bfl{{\bf l}}
\def\bfm{{\bf m}}
\def\bfn{{\bf n}}

\def\bfom{{\boldsymbol \omega}}

\def\bflam{{\boldsymbol \lambda}}
\def\bfeta{{\boldsymbol \eta}}
\def\bfxi{{\boldsymbol \xi}}

\def\ind{{\cal D}}
\def\nswv{{\cal C}_v}
\def\nsw{{\cal C}}
\def\Mmpq{M_m^{p,q}}
\def\phas{(x,\omega )}

\def\tii{{\cal T}_2}

\newcommand{\mif}{M^{\infty, 1 }}
\newcommand{\zd}{\Zst ^d}
\newcommand{\inv}{^{-1}}
\newcommand{\tfa}{time-frequency analysis}
\newcommand{\fif}{if and only if}
\newcommand{\rdd}{\Rst ^{2d}}
\newcommand{\riha}{Rihaczek}

\begin{document}

\title{\bf Pseudodifferential Operators on Locally Compact Abelian
  Groups and Sj\"ostrand's Symbol Class}

\author{Karlheinz Gr\"ochenig and Thomas Strohmer\footnote{K.~G.\ was
 supported by the Marie-Curie Excellence Grant MEXT-CT 2004-517154. 
T.~S.\ was supported by NSF DMS grants 0208568 and 0511461. $\qquad
\qquad \qquad \qquad \qquad \qquad \qquad \qquad $
  }}

\date{}
\maketitle







\begin{abstract}
We investigate pseudodifferential operators on  arbitrary locally compact abelian
groups. As symbol classes for the Kohn-Nirenberg calculus we introduce
a version of 
Sj\"ostrand's class. Pseudodifferential operators with
such symbols form a Banach algebra that is closed under
inversion. Since ``hard analysis'' techniques are not available on
locally compact abelian groups, a new time-frequency approach is used with
the emphasis on modulation spaces, Gabor frames, and  Banach algebras of
matrices. Sj\"ostrand's original results are thus understood as a
phenomenon of abstract harmonic analysis rather than ``hard
analysis'' and are  proved in their natural context and generality. 
\end{abstract}



\section{Introduction} \label{s:intro}

Pseudodifferential operators are a generalization of partial
differential operators, and the subject is usually treated with the
arsenal of ``hard analysis'', such as differentiation and
decomposition techniques, commutators etc. In this paper we develop a
new theory for pseudodifferential operators on general locally compact
groups instead of on $\rd $. Our main goal is to show the validity of
three subtle results of J.~Sj\"ostrand~\cite{Sjo94,Sjo95} on a Banach
algebra of pseudodifferential operators in the new context of locally
compact abelian groups. This is not a mere
generalization, because the formulation and extension of Sj\"ostrand's
results requires the development of  completely new methods in which
``hard analysis'' is replaced by phase-space (time-frequency)
analysis. 

To put the issues into a bigger context, recall  Wiener's Lemma: it
states   that a periodic function $f$ which has an absolutely summable  
Fourier series and which vanishes nowhere has an inverse $f^{-1}$ which 
also has an absolutely summable Fourier series~\cite{Wie32}.

This result can also be stated in the following way. Assume the sequence
$\{a_k\}_{k\in\Zst}$ satisfies $\sum_{k\in\Zst} |a_k| < \infty$. Let $A$ 
be a biinfinite Toeplitz matrix with entries $A_{k,l}=a_{k-l}$ for 
$k,l\in\Zst$. If $A$ is invertible on $\ltZ$ then its inverse $B:=A^{-1}$ 
has entries $ B_{k,l}=b_{k-l}$ which satisfy $\sum_{k\in\Zst} |b_k| < \infty$.
In this context $f(t)=\sum_{k\in\Zst} a_k e^{2\pi i kt}$ is called the 
{\em symbol} of $A$.

An intriguing generalization is due to Bochner and Philips~\cite{BP42} who
have shown that Wiener's Lemma remains true if the $a_k$ belong to a 
non-commutative Banach algebra instead of to $\Cst$.

In recent years several remarkable extensions of Wiener's Lemma have 
been published. Using results from~\cite{BP42},  Gohberg, Kaashoek, and 
Woerdeman~\cite{GKW89}, and independently Baskakov~\cite{Bas90} proved 
the following result. Consider the Banach algebra ${\cal C}$ of matrices 
$A$ with norm
\begin{equation}
\|A\|_{\nsw} := \sum_{k\in\\\Zst^d}\underset{i-j=k}{\sup}|A_{i,j}|
 < \infty.
\label{nswa}
\end{equation}
If $A \in \nsw$ is invertible in $\ltsp(\Zst^d)$ then its inverse 
$A^{-1}$ also belongs to $\nsw$.

Another Wiener-type theorem, this time in the context of pseudodifferential
operators, is due to Sj\"ostrand~\cite{Sjo95}. His striking result goes as 
follows. Let 
$g \in {\cal S}(\Rst^{2d})$ be a compactly supported $C^{\infty}$-function
satisfying the property $\sum_{k \in \Zst^{2d}} g(t-k) = 1$ for all
$t \in \Rst^{2d}$. Then a symbol $\sigma \in {\cal S}'(\Rst^{2d})$ belongs
to $M^{\infty,1}(\Rst^{2d})$ -- the {\em Sj\"ostrand class} -- if 
\begin{equation}
\int \limits_{\Rst^{2d}} \underset{k\in\Zst^{2d}}{\sup} 
|(\sigma\cdot g(. - k))^{\wedge}(\zeta)|\, d\zeta < \infty.
\label{sjostrandclass}
\end{equation}
Now let $K_{\sigma}$ be a pseudodifferential operator with (Weyl or
Kohn-Nirenberg) symbol $\sigma \in M^{\infty,1}(\Rst^{2d})$. (On $\rd
$, 
$K_\sigma $ is usually written as $\sigma (x,D)$ or $\sigma
(x,D)^w$). Sj\"ostrand 
proved that $M^{\infty ,1} $ is an algebra with respect to the
composition of pseudodifferential operators. Furthermore, if $\kn$ is
invertible on $\Ltsp(\Rst^{d})$, then  
$\kn^{-1} = K_{\tau}$ for some $\tau \in M^{\infty,1}(\Rst^{2d})$.
This is the Wiener property of $M^{\infty,1}(\Rst^{2d})$.
These results had a deep influence on recent work on new symbol
classes for pseudodifferential operators, as exemplified in the work
of Boulkhemair~\cite{boul99}, Lerner~\cite{ler05}, and Toft~\cite{Toft01}.

At first glance there is no relation between the two results on the
matrix algebra $\cC $ and the symbol class $M^{\infty , 1}$. However,  on
inspection of Sj\"ostrand's proof, which is in the realm of ``hard
analysis'', one sees that he uses the Wiener property of the
Gohberg-Baskakov matrix algebra $\cC $ as a tool. In fact, he found an
independent proof of their result, again using  a ``hard analysis'' approach with
commutators and decomposition methods.  On the other hand, 
 Gohberg et al.~and Baskakov use a ``natural'' approach in the context
 of harmonic analysis and  prove their result with classical
methods from Fourier analysis. 

One of the main insights of this paper is the observation that both
results are a manifestation  of a more general result, namely a Wiener
property  for a certain class of pseudodifferential operators on
locally compact abelian  (LCA)  groups. 

The main results of this paper can be summarized as follows:

(i) For every LCA group $\cG $ with  dual group  $\hat{\cG }$ we introduce a
symbol class $\mif (\tf )$. When $\cG = \rd $, this class reduces to
the Sj\"ostrand class, for $\cG = \zd $, this class coincides with the
matrix algebra $\cC$ defined in \eqref{nswa}. 

(ii) We show that $\mif (\tf
)$ is a Banach algebra under a twisted product that corresponds to the
composition of the corresponding pseudodifferential operators. 

(iii) We show that the Wiener property for the general Sj\"ostrand
class $\mif (\tf )$, i.e., if $\sigma \in \mif  (\tf )$ and
the pseudodifferential operator $K_\sigma $ is invertible  on $L^2
(\cG )$, then the inverse operator $K_\sigma \inv = K_\tau $ possesses
again a symbol $\tau \in \mif (\tf )$. 

(iv) We consider weighted versions of $\mif (\tf )$ and characterize
those weights for which the Wiener property holds.

The extension of Sj\"ostrand's original results to LCA groups is of
interest for both theoretical and practical reasons. \\
(a) Sj\"ostrand's  proof is based on
commutator estimates and decomposition techniques (``hard
analysis''). Such techniques  are not available on general LCA groups,
and it is by no means 
clear whether and how  such a generalization is actually possible. \\ 
(b) Proofs presented in a setting of LCA  groups show in some sense 
``what is really going on''.  The derivations are stripped off of 
lengthy analytic estimates and replaced by a time-frequency (phase
space) approach based on the ideas from~\cite{Gro04a}. Admittedly our 
approach requires more conceptual effort.  \\
(c) Sj\"ostrand's class and its weighted versions as well as the nonstationary
Wiener algebra have turned out to be very useful in applications, in
particular in the modeling of operators and transmission pulses
in connection with mobile communications, cf.~\cite{Str05}.
 The multidimensional setting is
potentially useful in applications such as spatially varying image or video 
(de)blurring. Furthermore, 
the  numerical implementation of pseudodifferential operators  requires a discrete 
finite setting, it is thus useful to know that the almost diagonalization
properties are preserved under appropriate discretization. 

(d) Pseudodifferential operators on the $p$-adic groups $\Qst _p$ occur
often in the construction of a $p$-adic quantum theory~\cite{RTVW89,Har93,Vla90}. Our
result hold in particular for operators on $\Qst_p$ and provide a new
type of a symbolic calculus. 

The paper is organized as follows. In Section~2 we develop 
time-frequency methods on locally compact abelian groups. This section
is somewhat lengthy, but we feel it
necessary to explain the main concepts of  time-frequency analysis,
such as amalgam spaces, modulation spaces, Gabor frames, and matrix
algebras. (By contrast,  in a paper
on standard pseudodifferential operators  it would  suffice  to refer to the expositions of
H\"ormander~\cite{hormander} or Stein~\cite{stein93}.)  In Section~3
we  explain the main formalism of pseudodifferential operators on
locally compact groups. Section~4 contains the key result about the
almost diagonalization of pseudodifferential operators in the
Sj\"ostrand class, and in Section~5 we formulate and prove our  main
results, the Banach algebra property and the Wiener property of
Sj\"ostrand's class on locally compact abelian groups. In  the final
Section~6 we discuss special groups and show that the matrix algebra
of Gohberg and Baskakov and Sj\"ostrand's symbol class are examples of
the same phenomenon. 

\section{ Tools from Time-Frequency Analysis}
\label{ss:tftools}

We first present the main concepts for \tfa\ on locally compact
abelian (LCA) groups. The constructions of \tfa\ are well-known for
$\rd $ and available in textbook form~\cite{Gro01,Fol89}. It is
less well known that \tfa\ works similarly on LCA groups. For some
contributions in this directions, we refer to \cite{FK98,Gro98}. 

In the following we focus on the details (weight functions on LCA
groups, spaces of test functions) that require special
attention. Whenever a  result can be formulated and proved as on $\rd
$, we will only formulate the result and refer to the proof on $\rd $.

\vspace{3 mm}

\noindent
\textbf{Locally Compact Abelian Groups.} Let $\group$ be a locally
compact abelian group. We assume that $\group$ 
satisfies the second countability axiom and is metrizable which is 
equivalent to the assumption that $\Ltsp(\group)$ is a separable Hilbert 
space~\cite{HR63}.
The elements of $\group$ will be denoted by italics  $x,y,u,\dots$ and
the group operation is written additively as $x+y$.
The dual group $\hatg$ is the set of characters on $\group$. We usually
denote characters by Greek letters $\xi, \eta , \omega,\dots$.
The action of a character $\xi  \in \hatg$ on an element $x\in\group$
is denoted by $\langle \xi,x \rangle$. Clearly the action of $-\xi$
on $x$ is then given by $\langle -\xi,x \rangle =
\overline{\langle \xi,x \rangle}$ where the overline denotes complex conjugation.

The phase-space or time-frequency plane is $\group \times \hatg$, its
elements are denoted by boldface letters $\bfx, \bfy, \bfu,\dots$.
By Pontrjagin's duality theorem~\cite{RS00} $\widehat{\hatg}$ is isomorphic
to $\group$, henceforth we will identify $\widehat{\hatg}$ with $\group$.
Consequently the dual group of $\tf$ is $\hattf$. We denote its 
elements by boldface Greek letters $\bfxi , \bfom ,\dots$.
Consistent with the previously introduced convention we
also write e.g.~$\bfx=(x,\xi) \in \tf $ and $\bfxi=(\xi,x) \in
\hattf$. 

By the structure theorem for locally compact abelian groups, 
$\group$ is isomorphic to a direct product $\group \simeq\Rst^d
\times \group_0$, where the LCA group $\group_0$ contains a compact open 
subgroup $\Ksub$~\cite{RS00,HR63}. Furthermore, if $\group_0$ contains the 
compact open subgroup $\Ksub$, then $\hatg_0$ contains the compact open 
subgroup $\Kperp$, cf.~Example~4.4.9 in~\cite{RS00} or 
Lemma~6.2.3 in~\cite{Gro98}.

The ``time-frequency plane'' (phase space) of $\group $ is $\group \times
\widehat{\group}$. As a consequence of  the structure theorem,
the phase-space is   $\group\times \hatg \simeq \Rst^{2d} \times
(\group_0 \times \hatg_0)$,  and $\group_0 \times \hatg_0$ contains
the compact-open group $\Ksub \times  \Kperp $.


The Fourier transform of a function $f$ on $\cG $ is 
is defined by~\cite{RS00}
\begin{equation}
\four f(\omega) = \hat{f}(\xi) = \int \limits_{\group} 
f(x) \overline{\langle \xi ,x\rangle} \, dx,\qquad 
\text{for $\xi \in \hatg$.}
\label{fourier}
\end{equation}
By Plancherel's Theorem $\cF $ is unitary from $L^2(\group )$ onto
$L^2 ( \hat{\cG } )$~\cite{RS00}. 

\vspace{3mm }

\noindent
\textbf{Time-Frequency Analysis.} 
For a function  $f$   on $\group$,  $x, y \in \group $, and $\xi \in
\hatg$,  we define the operators of 
{\em translation } $T_x$ and  {\em modulation } 
$M_{\xi}$ by
\begin{align}
T_y f(x)  = f(x-y), \qquad 
M_{\xi} f(x)  = \langle \xi ,x \rangle  f(x). 
\label{tmop}
\end{align}
The operators $T_x, M_{\xi}$ satisfy the commutation relations
\begin{equation}
T_x M_{\xi} = \overline{\langle \xi,x \rangle} M_{\xi} T_x.
\label{comm}
\end{equation}
The time-frequency shift operator $\pi(\bfx)$ on $\tf$ is defined by
$\pi(\bfx) = M_{\xi} T_{x}, \bfx=(x,\xi ) \in \tf$.

Given an appropriate function (``window'') $g$,  the 
{\em short-time Fourier transform} (STFT) of $f \in \Ltsp(\group)$ is 
defined by
\begin{equation}
\STFT_g f(x,\xi) = \int\limits_{\group} f(y) \overline{g(y-x)
\langle \xi,y \rangle}
dy,\qquad (x,\xi) \in \tf.
\label{stft}
\end{equation}
We note that
\begin{equation}
\STFT_g f(u,\omega) = \langle f, M_{\omega} T_u g \rangle  
= \langle \hat{f}, T_{\omega} M_{-u} \hat{g} \rangle = 
\STFT_{\hat{g}} \hat{f}(\omega,-u) \overline{\langle \omega,u \rangle}.
\label{STFTft}
\end{equation}
Furthermore, for $f, g \in L^2(\cG )$ there holds
\begin{equation}
\label{STFTfourier1}
\STFT_{M_{\eta} T_{y}g} M_{\xi} T_x f (u,\omega)=T_{(x-y,\xi-\eta)} 
\STFT_g f(u,\omega) \overline{\langle \omega - \xi, x \rangle} 
\langle \eta,u -x \rangle \, , 
\end{equation}
which follows from the definition of 
the STFT and the commutation relations~\eqref{comm}.

We will also make use of the following formula concerning the Fourier 
transform of a product of STFTs, which follows from an easy
computation (carried out  in~~\cite{rieffel88} and in~\cite{GZ01})
\begin{equation}
\label{eq:c1}
(V_{g_1}f_1{\overline{V_{g_2}f_2}})^{\wedge} \,
(\xi ,x )=(V_{f_2}f_1{\overline {V_{g_2}g_1}})(-x ,\xi ) \, .    
\end{equation}

\noindent
{\bf Weight Functions.} 
\begin{definition}
(a) A non-negative function $v$ on $\tf$ is called an {\em admissible weight}
if it satisfies the following properties:
\begin{itemize}
\item[(i)] $v$ is continuous, even in each coordinate,  and normalized such
that $v(0) = 1$.
\item[(ii)] $v$ is submultiplicative, i.e., $v(\bfx +\bfy)\le v(\bfx) v(\bfy)$,
$\bfx, \bfy \in \tf$.
\item[(iii)] $v$ satisfies the Gelfand-Raikov-Shilov (GRS) condition~\cite{GRS64}
\begin{equation}
\underset{n \toinf}{\lim} v(n\bfx)^{\frac{1}{n}} = 1 \qquad
\text{for all $\bfx \in \tf $}.
\label{GRS}
\end{equation}
\end{itemize}
(b) Let $v$ be an admissible weight. The class of $v$-moderate weights is 
\begin{equation}
{\cal M}_v = \Big\{m\ge 0: \underset{\bfx \in \tf}{\sup} 
\frac{m(\bfx+\bfy)}{m(\bfx)} \le Cv(\bfy), \quad \forall \,\bfy \in \tf\Big\}.
\label{moderateweights}
\end{equation}
\end{definition}

Examples: The standard weight functions on $\group $ are of the form
$$m(x) = e^{a \rho (x)^b} (1+\rho (x))^s \, , 
$$ where $\rho (x) = d(x,0)$
for some left-invariant metric $d$ on $\cG $.  Such a weight  is
submultiplicative, when  $a,s
\geq 0$ and $0\leq b \leq 1$, and $m$ satisfies the GRS-condition,
\fif\ $0\leq b < 1$. If $a, s \in \Rst $ are  arbitrary, then $m$ is $e^{|a|
  \rho (x)^b} (1+\rho (x))^{|s|}$-moderate. 

\vspace{3 mm}

\noindent
\textbf{Test Functions.}
For the treatment of weights of super-polynomial growth the standard space 
of test functions, the Schwartz-Bruhat space~\cite{Osborne}, is not suitable. We therefore 
introduce a {\em space of special test functions}, its construction is 
based on the structure theorem. Let $\Ksub$ be a compact-open subgroup of
$\group_0 $, let $\phi (x_1, x_2) = e^{-\pi x_1^2} \chiK(x_2) = 
\phi _1(x_1) \phi _2 (x_2)$ for $x= (x_1, x_2) \in \rd \times \cG _0$ and 
$$
\cS _{\cC } (\cG ) = \mathrm{span} \{ \pi (\mathbf{x})\phi :
\mathbf{x}\in \tf \} \subseteq L^2(\cG )  \,  
$$
be the linear space of all finite linear combinations of time-frequency
shifts of the ``Gaussian'' $\phi = \phi _1 \otimes \phi _2$. Then 
$$
V_\phi \phi (x,\xi ) =  V_{\phi _1} \phi _1 (x_1,\xi _1 ) \,   V_{\phi
  _2}\phi _2  (x_2 ,\xi _2 ) = e^{-\pi (x_1^2 + \xi _1^2)/2}  V_{\phi
  _2}\phi _2  (x_2 ,\xi _2 ) \, ,
$$ 
where $x_1, \xi _1 \in \rd $ and $(x_2, \xi _2) \in  \cG _0  \times \widehat{ \cG
_0}$.
Using the calculation on p.~228 of~\cite{Gro98} we find that 
\begin{eqnarray*}
   V_{\phi  _2}\phi _2  (x_2 ,\xi _2 ) &=& \langle \chiK, M_{\xi _2}
   T_{x_2} \chiK \rangle \\
&=& 
\begin{cases}
0 & \text{if}~x_2 \not \in \Ksub \\
( \chiK \cdot \chi _{x_2\Ksub})^{\wedge}\, (\xi _2)= 
\widehat{\chi}_\Ksub \, (\xi _2 ) = c(K) \chi _{\Kperp } (\xi _2) &
\text{if}~x_2 \in \Ksub\, .
\end{cases}
\end{eqnarray*}
Hence $V_{\phi _2} \phi _2 = c(\Ksub) \chiK \otimes \chi _{\Kperp } $,
where  $c(\Ksub)>0$ is a constant depending on $\Ksub$,   and  
the support of $ V_{\phi  _2}\phi _2$ is 
thus compact. Since a submultiplicative weight $v$ on $\rd $   grows at most
exponentially~\cite[Lemma VIII.1.4]{DS88}, we find that 
\begin{eqnarray}
  \label{eq:chr1}
\lefteqn{\int \int _{\tf } |V_\phi \phi (x,\xi )| v(x, \xi ) \, dx
  d\xi \leq } \\
&\leq & \int _{\rdd } e^{-\pi (x_1^2 + \xi _1^2)/2} \, v(x_1, 0 , \xi
_1, 0) \, dx_1
d\xi _1 \, \, \int _{\cG _0 \times \widehat{\cG _0}} \chiK (x_2) 
\chi_{\Kperp} (\xi _2) \, v(0,x_2, 0, \xi _2) \, dx_2 d\xi _2 < \infty \,
  . \notag 
\end{eqnarray}
 Consequently, $V_\phi \phi $ and hence every function in 
$\cS _{\cC} $ is integrable with respect to arbitrary moderate weight
functions. 

\vspace{3mm}

\noindent {\bf Modulation Spaces.} Let $m$ be a weight function.
We define the $\Mmpq $-norm of $f \in \cS _{\cC } (\cG )$ to be 
\begin{equation}
  \label{eq:chr6}
  \|f\|_{\Mmpq}:=\|V_{\varphi } f \, m \|_{L^{p,q}} = \biggl( \int
  _{\hatg } \biggl(
    \int _{\cG } |V_{\varphi } f (x,\xi )|^p\, 
m(x,\xi)^p \, dx\biggr)^{q/p} d\xi\biggr)^{1/q}\, .
 \end{equation}

Analogous to~\cite{Gro04a} we define the {\em modulation space} $\Mmpq
(\cG )$ 
as the completion of the space $\cS _{\cC} (\cG )$ with respect to the 
$\Mmpq $-norm, when $p,q < \infty $ and the weak$^*$-completion if 
$pq=\infty$. If $p=q$ we also write $M^{p}_m$ instead of $M^{p,p}_m$.
If $m\geq 1$, then $\Mmpq $ is a subspace of $M^{\infty} $ which 
in turn is a particular subspace of the space of tempered distributions 
$\cS '$ on $\cG $. 

In the sequel we will distinguish between modulation spaces on $\cG $
and on $\tf $. The space  $\modiiv(\group)$ will serve as 
``window space''  and can be considered a space of test
functions. Pseudodifferential operators will act on the modulation
spaces $M^{p,q}_M(\cG )$.

Modulation spaces on $\tf$ will be used
as symbol classes. In particular,  $\modv(\tf)$ is the appropriate
generalization of Sj\"ostrand's class to LCA groups.
If we take the norm completion of 
$\cS _{\cC}$ in the $M^{\infty, 1}_v$-norm, we obtain a class of symbols 
that leads to compact operators, see~\cite{BGH05}.

We will use the following standard properties of modulation spaces. 

\begin{proposition}[Duality]
  \label{duality}
(i) Let  $1\leq p,q< \infty $ and $p' = \frac{p}{p-1}$ be the conjugate
index. Then the dual space of $M^{p,q}_m(\cG ) $ is the modulation
space $M^{p',q'}_{1/m} (\cG )$.

(ii)  $M^{1}_m $ is the dual space of
$M^{0,0}_{1/m} := \mathrm{clos}_{M^{\infty, \infty}_{1/m} }(\cS _{\cC })$ 
(the closure of the test functions with respect to the $M^{\infty} $-norm), likewise $M^{\infty
  ,1}_m (\cG )$ is the dual of $\mathrm{clos}_{M^{1,\infty }_{1/m}}(\cS _{\cC
})$ (See~\cite{BGH05}).

(iii)  If $f \in M^{\infty ,1}_m (\cG )$ and $g \in M^{1,\infty
}_{1/m}$, then $\langle f,g\rangle := \int _{\tf } V_\varphi f (\bfx )
\, \overline{V_\varphi g(\bfx )} \, d\bfx $ is well-defined and
satisfies
\begin{equation}
  \label{eq:chr5}
  |\langle f, g \rangle | \leq C\|f\|_{M^{\infty ,1}_m (\cG )} \, \|g\|_{ M^{1,\infty
}_{1/m}}\, .
\end{equation}
\end{proposition}

\begin{lemma}\label{equivnorm}
  If $g\in M^{1}_v$ and $m \in {\cal M} _v $, then $\|V_g f \, m
  \|_{L^{p,q}}$ is an equivalent norm on $\Mmpq $~\cite[Ch.~11]{Gro01}. 
\end{lemma}
%


\medskip
\noindent
{\bf Amalgam Spaces.} 
A \emph{lattice} $\Lambda $ of $\cG $ is a discrete subgroup such that $\cG
/\Lambda $ is compact. Then there exists a relatively compact set
$U\subseteq \cG $, a \emph{fundamental domain} for $\Lambda $, such
that $\bigcup _{\lambda \in \Lambda } (\lambda
+U)  = \cG $ and $(\lambda + U) \cap (\mu +U) = \emptyset$ for
$\lambda \neq \mu \in \Lambda $.

If $\cG $ does not have a lattice, as  is  the case for $p$-adic
groups, we resort to the following construction. Recall  the structure
theorem $\cG \simeq  \Rst ^d \times \cG _0$, where $\cG _0 $ possesses
the compact-open subgroup $\Ksub$. Now choose an invertible, real-valued
$d\times d$-matrix $A$, and a set of coset representatives 
$D$ of $\group_0/\Ksub$ in
$\cG _0$, and let $U= A[0,1)^d \times \Ksub$.  Then $\cG  = \bigcup
_{\lambda \in \Lambda } (\lambda +U)$ is a partition of $\cG $. We
call the discrete set $\Lambda := A\Zst ^d \times D$ a \emph{quasi-lattice} with \emph{fundamental domain} $U$. 

Consequently a quasi-lattice in the time-frequency plane $\tf $ 
will have the form $\Lambda = \Lambda _1 \times \Lambda _2 := (A\rd
\times D_1) \times (B\rd \times D_2) \simeq \cA \rdd \times D_1 \times
D_2$ (where $A, B$ are $d\times  d$
invertible matrices) with fundamental domain $U=U_1\times U_2 =
(A[0,1)^d \times \Ksub) \times (B[0,1)^d \times \Ksub ^\perp )$. 

Using this construction, we can now define amalgam spaces on $\tf $,
see~\cite{FS85} and \cite{feichtinger90} for a detailed theory. 
 
\begin{definition}
\label{amalgam}
Let $\Lambda$ be a quasi-lattice of $\tf $ and $U$ a relatively compact
fundamental domain of $\Lambda$ in $\tf $. Let $m$ be a weight
function on $\tf $.
A continuous function $F$ on $\tf$ belongs to the {\em amalgam space} 
$W(C,\ell^{p,q}_m)(\tf )$ if the sequence $\{a(\bfl )\}_{\bfl  \in \Lambda}$ with
\begin{equation}
a(\mathbf{l} ) = a(l, \lambda)  = \underset{(u,\eta ) \in
  U}{\sup} |F(u+l, \eta  + \lambda  )|
\end{equation}
belongs to $\ell^{p,q}_m(\Lambda)$, that is 
$\Big( \sum _{\lambda \in \Lambda _2} \big( \sum _{l\in \Lambda _1}
a(l,\lambda )^p m(l,\lambda )^p \big)^{q/p} \Big)^{1/q}< \infty $, with the  usual
modifications if $p \, q=\infty$. 
\end{definition}
We note that the definition 
of the amalgam spaces is independent of the quasi-lattice $\Lambda$ and the
fundamental domain $U$, and 
different choices for $\Lambda$ lead to equivalent norms~\cite{FS85}.

Among others, amalgam spaces occur in \tfa\ in the description of the
fine local properties of the STFT. 

\begin{theorem}\label{wcamalg}
  Assume that $g\in M^1_v (\cG )$, $f\in \Mmpq (\cG)$ for $1\leq p,q
  \leq  \infty  $, and $m$ 
  a $v$-moderate weight. Then $V_g f $ is in $W(C,\ell ^{p,q}_m )(\tf )$. In
  particular, $V_gg \in W(C, \ell ^1_v ) (\tf )$. 
\end{theorem}

\begin{proof}
  The statement is a special case of~\cite[Lemma~7.2, Thm.~8.1]{fg89mh}
  (use the representation $(x,\xi ,\tau ) \to \tau T_x
  M_\xi $ on $L^2(\cG )$ of the Heisenberg-type group $\tf \times \mathbb{T}
  $). A direct proof for $\rd $ is given in \cite[Thm.~12.2.1]{Gro01}. 
\end{proof}

As an important consequence, we formulate this result for the general
Sj\"ostrand 
class $M^{\infty ,1}_v(\tf )$, where $v$ is an admissible weight on
$\hattf$. Note that $V_\Phi \sigma $ is a function on $( \tf ) \times
(\hattf )$.

\begin{corollary}\label{amalcor}
 Let $\tilde{\Lambda }$ be a quasi-lattice in $\hattf $
with fundamental domain $\tilde{U}$. If $\Phi \in M^1_{1\otimes v} (\tf )$ and
$\sigma \in M^{\infty ,1} _{v} (\tf )$, then the sequence 
$$
h(\bflam ) := \sup _{\bfeta \in \tilde{U} } \sup _{\bfx \in \tf }
|V_\Phi \sigma (\bfx ,  \bflam + \bfeta  )|
$$
is in $\ell ^1 _v (\tilde{\Lambda })$. 
\end{corollary}

\medskip
\noindent
{\bf Gabor Frames.}
We assume familiarity with Gabor frames and refer to~\cite{Gro01},
Ch.~5 and 7, for details.
Given a quasi-lattice $\Lambda \subset \tf$ and a {\em window} $g \in
\Ltsp(\group)$  
the associated {\em Gabor system} $\{g_{m,\mu}\}_{(m,\mu)\in\Lambda}$
consists of functions of the form
\begin{equation}
g_{m,\mu} = M_{\mu} T_{m} g,\,\,\, (m,\mu) \in \Lambda.
\label{gaborsystem}
\end{equation}
The analysis operator or coefficient operator
$C_g: \Ltsp(\group) \mapsto \ell^2(\Lambda)$ is defined as
\begin{equation}
\label{analysis}
C_g f = \{\langle f, M_{\mu} T_{m} g \rangle\}_{(m,\mu) \in \Lambda}.
\end{equation}
The adjoint operator, which is also known as synthesis operator,
can be expressed as
\begin{equation}
C_g^{\ast} \{c_{m,\mu}\}_{(m,\mu) \in \Lambda} = 
\sum_{(m,\mu) \in \Lambda} c_{m,\mu} M_{\mu} T_{m} g \qquad 
\text{for $\{c_{m,\mu}\}_{(m,\mu) \in \Lambda}\in\ell^2(\Lambda)$.}
\label{synthesis}
\end{equation}
Associated to a Gabor system is the Gabor frame operator $S$  defined as 
\begin{equation}
Sf = \sum_{(m,\mu) \in \Lambda} \langle f,M_{\mu}T_{m} g \rangle
M_{\mu}T_{m} g = C_g ^* C_g f \, .
\label{frameop}
\end{equation}
We say that $\gabframe$ with $g\in \Ltsp(\group)$ is a Gabor frame for 
$\Ltsp(\group)$ if $S$ is invertible on $L^2 (\cG )$. Equivalently
there exist constants $A,B>0$ such that 
\begin{equation}
A \|f\|_2^2 \le \sum_{(m,\mu)\in\Lambda}|\langle f,M_{\mu}T_{m}g 
\rangle|^2 = \langle Sf,f \rangle \le B \|f\|^2_2, \qquad \text{for all $f \in \Ltsp(\group)$.} 
\label{framedef}
\end{equation}
A Gabor system $\gabframe$ is called a {\em tight} Gabor frame if $A=B$ in~\eqref{framedef}.
In this case $S$ is just (a multiple of) the identity operator on
$\Ltsp(\group)$.

 For our purposes we need tight Gabor frames
generated by a window $g \in \modiivg$. 
The existence and  construction of Gabor frames are well understood on
$\rd $, but our knowledge of explicit Gabor frames on LCA groups is
thin. Therefore the  following existence theorem may be of independent
interest.  

\begin{theorem}\label{gabexist}
  Let $v$ be an admissible weight on $\tf  $ satisfying the
  GRS-condition, and let $\Lambda := \alpha \mathrm{I} \times D
  $ be a quasi-lattice in  $ \tf $ with $\alpha <1$ and $D$ a set
  of representatives of $\cG/ \Ksub \times \widehat{\cG _0} / \Ksub
  ^\perp $. Then there   exists a $g \in M^1_v (\cG )$, such that $\{
  \pi (\lambda )g: \lambda 
  \in \Lambda \}$ is a tight Gabor frame for $L^2(\cG )$. 
\end{theorem}

\begin{proof}
  According to the structure theorem we distinguish several cases. 

\textbf{Case I: $\cG = \rd $.} We choose $\Lambda = \alpha \zdd $ for
$\alpha <1$ and the Gaussian  $\phi (t)= e^{-\pi t\cdot t}$.
It follows from the main result in~\cite{Lyu92,SW92}
that $\{M_{\mu}T_m \phi \}_{(m,\mu) \in \Lambda}$ 
is  a Gabor frame for $L^2(\rd)$ with $\phi \in \modiivg$. 
To this Gabor frame we apply Cor.~4.5 of \cite{GL02}: \emph{Let
  $v$ be an admissible weight satisfying the 
GRS-condition. If   $\gabframe$ is a Gabor frame for $L^2(\Rst^d) $ with  
$g\in \modiiv(\Rst^{d})$ and associated frame operator $S$,  then $\{
M_\mu T_m S^{-1/2}g: (m,\mu ) \in \Lambda \} $ is a  tight frame
and the window $\gamma ^\circ = S^{-1/2}g $ 
also belongs to $\modiivg$.} This construction provides an abundance of
tight Gabor frames for $L^2(\rd )$. 

\textbf{Case II: $\cG = \cG _0$}, where $\group _0$ contains the
compact-open subgroup $\Ksub $.  Let $D_1$ be a set of coset
representatives of $\group_0/\Ksub$ and  $D_2$ be a set of coset
representatives of $ \hatg_0/\Kperp$. Then 
 $D = D_1\times D_2$ is a quasi-lattice in $\tf $,  and  the family
 $\{M_{\delta} T_d \chiK :(d, \delta) \in D\}$  
is an orthonormal basis for $L^2(\group_0)$.   To verify this claim,
we note that $\{\delta : \delta \in \widehat{\Ksub}  \}$ is an
orthonormal basis for $L^2(\Ksub )$, because $\Ksub $ is
compact. Since $\widehat{\Ksub } \simeq \widehat{\group _0}/ \Ksub
^{\perp} \simeq D_2$, the set $\{ M_\delta \chi _{\Ksub }: \delta \in
D_2\}$  is an orthonormal basis for $L^2(\Ksub) \subseteq L^2 (\cG
_0)$. Furthermore, since the translates $T_{d_1} \chiK, T_{d_2} \chiK$
have disjoint support for $d_1 , d_2 \in D_1 \simeq \hatg_0/\Ksub,
d_1\neq d_2$ and since 
$$L^2(\group_0) = \oplus_{d \in D_1} L^2(d\, \Ksub),$$
 $\{M_{\omega}T_d \chiK\}_{(d,\delta) \in D}$ is an 
orthonormal basis for $L^2(\group_0)$.

Furthermore, $\chiK\in \modiiv(\group_0)$ by Proposition~6.4.5 
in~\cite{Gro98} or as a consequence of~\eqref{eq:chr1}. 

\textbf{Case III: $\cG \simeq \rd \times \cG _0$ is an  arbitrary LCA
  group.} Let $\Lambda = \alpha \zdd \times D$, $\alpha < 1$ be a
quasi-lattice in $\tf $, let $\{M_\mu T_m \gamma ^\circ : (m,\mu
) \in \alpha \zdd \}$ be a tight frame for $L^2(\rd )$, and $
\{M_{\delta} T_{d} \chi _{\Ksub } :  (d,\delta
) \in  D \}$ be the orthonormal basis for $L^2(\cG _0)$. Then the
set $\{ M_{(\omega,\delta )} T_{(u,d )} (\gamma ^{\circ} \otimes
\chi _{\Ksub }) \}$ is a tight frame for $L^2(\group )$, because the 
tensor product of (tight) frames is again a (tight) frame. Finally,
the window $g = \gamma ^{\circ } \otimes \chi _{\Ksub }$ is in $M^1_v
(\cG )$, which is shown as in \eqref{eq:chr1}.
\end{proof}

\medskip
\noindent
{\bf A  Banach Algebra of Matrices.} 
The following matrix algebra is a natural generalization of Wiener's
algebra and  will play a central role in our investigations.
\begin{definition}[\cite{GKW89,Bas90,Bas97}]
\label{def:nsw}
 Let $\ind $ be   a countable discrete abelian  subgroup, and let $v$ be an
admissible weight on $\ind $.  
The nonstationary Wiener algebra $ \nswv = \nswv(\ind)$ consists of all matrices 
 $A = [A_{i,j}]_{i,j \in \ind \times \ind}$ on the index set $\ind$,  for which
\begin{equation}
\|A\|_{\nswv(\ind)} := \sum_{j\in\ind}\,
\underset{i\in \ind }{\sup}|A_{i,i-j}|\, v(j)
\label{nswv1}
\end{equation}
is finite.
\end{definition}
It is easy to verify that
\begin{equation}
\sum_{j\in\ind}\underset{i}{\sup}|A_{i,i-j}|v(j) =
\underset{a \in \ell^1_v(\ind)}{\inf}\{\|a\|_{\ell ^1_v} : |A_{i,j}| \le a(i-j), i,j\in \ind\}.
\label{nswv2}
\end{equation}
The unweighted version of the following result was mentioned
already in the introduction of this paper.  
\begin{theorem}
\label{th:bask}
Let $A = [A_{i,j}]_{i,j \in \ind}$ be a matrix in 
$\nswv$ where $v$ is an admissible weight. If $A$ is invertible
on $\ell^2(\ind)$,  then $A^{-1} \in \nswv$.
\end{theorem}
This theorem was obtained by Gohberg et al.~\cite{GKW89} and 
independently by Sj\"ostrand~\cite{Sjo95} for the case $v=1$. 
The weighted case as well as quantitative versions 
were derived by Baskakov~\cite{Bas90,Bas97}. 

We need the following generalization of this theorem. 
We recall that a matrix $A $ possesses a pseudoinverse, if there
exists a subspace $\mathcal{M} \subseteq \ell ^2 (\ind )$, such that
the restriction of $A$ to $\mathcal{M}$ is invertible and
$\mathrm{ker} (A) = \mathcal{M}^\perp $. The trivial extension of this
inverse on $\mathcal{M}$ to all of $\ell ^2 (\ind )$ is called the
pseudoinverse and denoted by $A^{+}$. 
\begin{corollary}[\cite{FG04}]
\label{cor:pinv}
Let $A = [A_{i,j}]_{i,j \in \ind}$ be a matrix in 
$\nswv$ where $v$ is an admissible weight. If $A$ has a pseudoinverse
$A^{+}$, then $A^{+} \in \nswv$.
\end{corollary}

\section{Pseudodifferential Operators on Locally Compact Abelian
  Groups} \label{s:lca}

We turn to the investigation of pseudodifferential operators on LCA
groups. The abstract formalism was developed in~\cite{FK98}. With
proper notation, most formulas are almost identical to those for
pseudodifferential operators on $\rd $. 
\begin{definition}
Let $\sigma$ be a function or distribution in $M^{\infty} (\tf )$. 
The pseudodifferential 
operator with Kohn-Nirenberg symbol $\sigma$ is the operator $\kn$ given by
\begin{equation}
(\kn f)(x) = \int \limits_{\hatg} \sigma(x,\xi)
\hat{f}(\xi) \langle \xi,x \rangle d\xi.
\label{kohn}
\end{equation}
If $F$  is some function space, we write $\kn \in \Op(F)$ whenever 
$\sigma \in F$.
\end{definition}
Alternatively we can write $\kn $ as a superposition of time-frequency
shifts~\cite{FK98,Fol89,Gro01}: 
\begin{equation}
\label{spreading}
\kn f(x)=\int\limits_{\hattf}\hat{\sigma}(\omega,u) M_{\omega}
T_{-u}f(x)d\omega du \, .
\end{equation}
If $ \hat{\sigma } \in L^1 (\hat{\cG } \times \cG )$, $f \in
L^1(\cG )$, and $\hat{f} \in L^1 (\hatg )$, 
 this follows from the computation
\begin{gather}
\kn f(x)=\int \limits_{\hatg} \sigma(x,\xi)
\hat{f}(\xi)  \langle\xi,x \rangle d\xi  
=\int \limits_{\tf} \sigma(x,\xi)f(y) \overline{\langle\xi,y-x\rangle} 
dy d\xi \notag \\ 
=\int\limits_{\hattf}\hat{\sigma}(\omega,y-x) f(y)\langle\omega,x\rangle
d\omega dy
=\int\limits_{\hattf}\hat{\sigma}(\omega,u) f(x+u)\langle\omega,x \rangle
d\omega du \notag \\
=\int\limits_{\hattf}\hat{\sigma}(\omega,u) M_{\omega} T_{-u}f(x) d\omega du.
\end{gather}
Expression~\eqref{spreading} is called the {\em spreading 
representation} of $\kn$ and $\hat{\sigma}$ is the {\em spreading function}.
For more general symbol classes the validity of the spreading
representation follows by a routine density argument~\cite{FK98,Gro04b}.
Expression~\eqref{spreading} represents pseudodifferential
operators as linear combination of time-frequency shift operators,
which suggests that methods from time-frequency analysis are a natural
tool for the study of pseudodifferential operators.

We have the following formal symbol calculus for Kohn-Nirenberg 
pseudodifferential operators. 

\begin{lemma} 
\label{twistedlemma}
If $ \hat{\sigma }, \hat{\tau } \in L^1 (\hat{\cG } \times \cG )$, then 
\begin{equation}
\kn \knt = \knst,
\label{symbolcalculus}
\end{equation}
where the twisted convolution $\twist$ of $\hat{\sigma},\hat{\tau}$ is
defined by 
\begin{equation}
\hat{\sigma} \twist \hat{\tau}(\xi,u) = \int \limits_{\hattf}
\hat{\sigma}(\zeta,y)\hat{\tau}(\xi-\zeta,u-y) \langle\xi-\zeta,y \rangle
d\zeta dy.
\label{twisted}
\end{equation}
\end{lemma}

\begin{proof}
Our hypothesis guarantees that the integrals below converges
absolutely and thus Fubini's theorem permits to change the order of
integration. 
\begin{align*}
\label{}
\kn \knt f &= \int \limits_{\hattf}\hat{\sigma}(\zeta,y)M_{\zeta}T_{-y}
d\zeta dy 
\int \limits_{\hattf} \hat{\tau}(\xi,u) M_{\xi}T_{-u} f d\xi
du \\
&= \int \limits_{\hattf} \int \limits_{\hattf}
\hat{\sigma}(\zeta,y) \hat{\tau}(\xi,u) \, \langle \xi,y \rangle
\, M_{\zeta+\xi}T_{-(y+u)} 
f d\xi du d\zeta dy \\
&= \int \limits_{\hattf} \Big( \int \limits_{\hattf}
\hat{\sigma}(\zeta,y) \hat{\tau} (\xi-\zeta,u-y )  \langle\xi-\zeta, y\rangle 
d\zeta dy \Big) \, M_{\xi}T_{-u} f d\xi du  \\
&= \int \limits_{\hattf} (\hat{\sigma}\twist\hat{\tau})(\xi,u)
M_{\xi}T_{-u} f d\xi du  = \knst f.
\end{align*}
\end{proof}

\begin{definition}
Let $v$ be an admissible weight on $\hat{\cG } \times \cG $. 
The weighted Sj\"ostrand class $\Op(\modv(\tf))$ is the class of 
pseudodifferential operators $\kn$ whose symbol 
$\sigma \in M^\infty  (\tf )$ satisfies
\begin{equation}
\|\sigma\|_{\modv} =  \int \limits_{\hattf}
\underset{\bfx \in \tf}{\sup} |\STFT_{\Psi} \sigma(\bfx,\bfom)|\, v(\bfom)
d\bfom < \infty
\label{modulinf1}
\end{equation}
with $\Psi \in  {\cal S}_{\mathcal{C} }(\tf)  \setminus \{0\} $.
\end{definition}

Note that the weight in~\eqref{modulinf1} depends only on $\bfom $, so 
consistency with  \eqref{eq:chr6} would require the clumsier notation
$M^{\infty , 1} _{1\otimes v} (\tf )$. 

Let $\gabframe$ be a Gabor system for $\Ltsp(\group)$ with respect to
a quasi-lattice $\Lambda \subseteq \tf$. Using the notation
of time-frequency shift operators $\pi (\bfx ) = M_\xi T_x$,  we can
write this system as 
$\{\pi(\bfm) g\}_{\bfm \in \Lambda}$. For a given
pseudodifferential operator $\kn$ we define the matrix $M(\sigma)$ by
\begin{equation}
[M(\sigma)]_{\bfm,\bfn} = \langle \kn \pi(\bfn )g,\pi(\bfm )g \rangle,
\qquad \bfm,\bfn \in \Lambda.
\label{matrix1}
\end{equation}
Since $M(\sigma)$ depends also on $g$ and $\Lambda$, it would be more
accurate to use the notation $M(\sigma,g,\Lambda)$. However, whenever
there is no danger of confusion,  we will simply write $M(\sigma)$.

Assume $\kn$ is bounded on $\Ltsp(\group)$ and that $\piframe$ is a tight 
frame for $\Ltsp(\group)$ with (lower and upper) frame bound equal to 1. 
In this case $C^{\ast} C = I$, where $C_g$ and $C_{g}^{\ast}$ are defined 
as in~\eqref{analysis} and~\eqref{synthesis}. We can represent 
$f\in L^2(\group)$ as 
$f = \sum_{\bfn \in \lambda} \langle f,\pi(\bfn) g\rangle \pi(\bfn)g$.
For $\bfm \in \Lambda$ we compute
\begin{align}
C_g (\kn f)(\bfm) & = \langle \kn f,\pi(\bfm) g \rangle \notag \\
& = \sum_{\bfn \in \Lambda} \langle f,\pi(\bfn) g \rangle
\langle \kn \pi(\bfn) g,\pi(\bfm) g \rangle = (M(\sigma) C_g f)(\bfm).
\label{knmat}
\end{align}
Since $C^{\ast} C = I$, equation~\eqref{knmat} can be expressed equivalently as
\begin{equation}
\kn f = C_{g}^{\ast} M(\sigma)C_g f.
\label{matrix2}
\end{equation}

The following lemma specifies the kernel and the range of $M(\sigma )$
and is taken from~\cite[Lemma~3.4]{Gro04a} (the  proof carries over almost verbatim 
to our setting  by replacing the Weyl symbol by the Kohn-Nirenberg
symbol and $\Rst^{2d}$ by $\tf$).
\begin{lemma}
\label{le:range}
Let $\gabframe$ be a Gabor frame for $\Ltsp(\group)$. If $\kn$ is bounded on 
$\Ltsp(\group)$ then 
$M(\sigma)$ is bounded on $\llsp$ and maps $\ran(C_g)$ into $\ran(C_g)$ 
with $ \ran(C_g)^{\perp} = \ker(C^{\ast}_g) \subseteq \ker(M(\sigma))$. 
\end{lemma}

\section{Almost Diagonalization} 

We characterize the Sj\"ostrand class by its almost diagonalization
with respect to Gabor frames. The corresponding results on $\rd $ were
obtained in~\cite{Gro04a} and in a slightly different version that is
more suitable to applications in~\cite{Str05}.

\begin{definition}
The (cross) Rihaczek distribution of $f,g \in \Ltsp(\group)$ is defined as
\begin{equation}
R(f,g)(x,\xi )  = f(x) \overline{\hat{g}(\xi )} \,\, 
    \overline{\langle \xi, x\rangle}.
\label{rihaczek}
\end{equation}
\end{definition}

The next lemma states two important properties of the Rihaczek
distribution and clarifies its appearance  in the analysis of the
Kohn-Nirenberg pseudodifferential operators.  Now define    
\begin{equation}
  \label{eq:chr3}
  \cJ (x,\xi ) =  (-\xi , x ) \qquad (x,\xi ) \in \tf \, ,
\end{equation}
then $\cJ $ is an  isomorphism from $ \tf $  onto $\hattf$, and $\cJ
$ preserves the Haar measure.  

\begin{lemma}
\label{key}
Let $f, g \in \cS_{\cC}(\group)$.  Then:\\
(i) 
\begin{equation}
\label{riha1}
  R(\pi (\mathbf{x})g, \pi (\mathbf{y})f) = \langle \eta , x-y\rangle \,
  M_{\cJ (\mathbf{y}- \mathbf{x})}   T_{(x,\eta )} R(g,f),
\end{equation}
(ii)
\begin{equation}
\langle \kn  \pi(\bfx)f, \pi(\bfy) f \rangle = 
\overline{\langle \eta,x-y \rangle } \STFT_{R(g,f)} \sigma
\Big((x,\eta),{\cal J}(\bfy-\bfx)\Big). 
\label{rihaSTFT}
\end{equation}
(iii) If $f,g \in M^1_v (\group)$ for  an admissible weight $v$  on
$\tf $, then \\ $R(g,f) \in 
M^1 _{1\otimes v\circ \cJ \inv }(\tf )$. 
\end{lemma}

\begin{proof}
%
(i) follows from the calculation 
\begin{eqnarray*}
  R(\pi (\mathbf{x})g, \pi (\mathbf{y})f)(t,\tau ) &=& \pi (\mathbf{x})
g(t) \,
  \overline{(\pi (\mathbf{y})f )\, \widehat{}\, (\tau )}\,\,
  \overline{\langle\tau, t\rangle}  \\
&=& M_\xi T_x g(t) \, \overline{T_\eta M_{-y} \hat{f} (\tau )}\,
\overline{\langle\tau, t\rangle}  \\
&=& \langle \xi , t\rangle g(t-x) \langle  \tau - \eta , y \rangle
\overline{\hat{f} (\tau - \eta )}  \\
& & \quad \cdot \overline{\langle \eta , t\rangle}
\overline{\langle \tau , x\rangle } \langle \eta , x\rangle
\overline{\langle \tau - \eta , t-x \rangle} \\
&=& \langle \eta , x-y\rangle \, M_{(\xi - \eta , y-x)} T_{(x,\eta )
}\,R(g,f)(t,\tau ) \, .
\end{eqnarray*}
Using the definition for $\cJ $, we have obtained \eqref{riha1}


(ii) We first calculate the action of $K_{\sigma}$ for functions 
$f,g\in\cS_{\cC}(\group)$ in terms of the Rihaczek distribution:
\begin{eqnarray}
  \langle  K_\sigma f, g \rangle &=& \int \limits_{\group} \int \limits_{\hatg}
  \sigma (x,\xi ) \hat{f}(\xi ) \langle \xi , x\rangle \,
  \overline{g(x)}\, d\xi dx  \notag \\
&=& \langle \sigma , R(g,f) \rangle.
\label{kohnriha1}
\end{eqnarray}
Using~\eqref{kohnriha1} and~\eqref{riha1} we compute
\begin{eqnarray}
  \langle K_\sigma \pi (\mathbf{y})f, \pi (\mathbf{x})g \rangle &=&
  \langle \sigma ,  R(\pi (\mathbf{x})g, \pi (\mathbf{y})f))\rangle \notag \\
&=& \overline{\langle \eta , x-y\rangle }\, \langle \sigma , M_{\cJ (\mathbf{y}-
\mathbf{x})}
  T_{(x,\eta )} R(g,f) \rangle \notag\\
&=&\overline{ \langle \eta , x-y\rangle } \,  V_{R(g,f)} \sigma ((x,\eta ), \cJ
(\mathbf{y} - \mathbf{x})) \, .
\label{kohnriha2}
\end{eqnarray}
which is~\eqref{rihaSTFT}.

(iii) is proved  similar to Lemma~\ref{wdstft}(ii) and therefore omitted. 
\end{proof}

Now we state the key theorem about almost diagonalization of
pseudodifferential operators $\kn$ for symbols in the generalized
Sj\"ostrand class $ \modvj(\tf)$. First we formulate a version for LCA
groups that contain a lattice $\Lambda $ and a corresponding Gabor
frame.  
\begin{theorem}
\label{th:approxdiag}
Let $g \in \modiiv(\group)$ for some  admissible weight $v$ on $\tf $, $\Lambda
\subseteq \tf $ be a lattice, 
and assume that $\piframe$ is a tight Gabor 
frame for $\Ltsp(\group)$. 
Then for $\sigma \in  M^{\infty} (\tf )$ 
the  following properties are equivalent.
\begin{itemize}
\item[(i)] $\sigma \in \mod_{v\circ \cJ^{-1}}(\tf)$.
\item[(ii)] There exists a function $H \in \Lsp^1_v(\tf)$ such that
\begin{equation}
|\langle \kn \pi(\bfz )g,\pi(\bfy)g \rangle| \le H(\bfy - \bfz ),
\qquad \text{for all $\bfy,\bfz \in \tf$.}
\label{}
\end{equation}
\item[(iii)]
There exists a sequence $h \in \lsp^1_v( \Lambda)$ such that
\begin{equation}
|\langle \kn \pi(\bfn)g ,\pi(\bfm )g \rangle| \le h(\bfm-\bfn),
\qquad \text{for all $\bfm,\bfn \in \Lambda$.}
\label{}
\end{equation}
\end{itemize}
\end{theorem}
\begin{proof} \\
$(i) \Rightarrow (ii)$.  Let $\sigma \in \mod_{v\circ \cJ^{-1}}(\tf)$. 
Denote $\bfx = (x,\xi)$ and $\bfy = (y,\eta)$, and set $\Psi=R(g,g)$,
which is in $M^1_{1\otimes v\circ \cJ^{-1}}(\tf )$ by
Lemma~\ref{key}(iii). 
We use Lemma~\ref{key}(ii) to compute
\begin{equation}
|\langle \kn  \pi(\bfx)g, \pi(\bfy) g \rangle|  = 
|\STFT_{\Psi} \sigma \big( (x,\eta), {\cal J}(\bfy-\bfx) \big)| \le 
\underset{\bfz \in\tf}{\sup}|\STFT_{\Psi}\sigma(\bfz,{\cal J}(\bfy-\bfx))|.
\label{eststft}
\end{equation}
Now set $H (\bfx):=\underset{\bfz \in\tf}{\sup} |\STFT_{\Psi}\sigma(\bfz,\mathcal{J}\bfx)|$.
Then 
\begin{eqnarray*}
  \int _{\tf } H(\bfx ) v(\bfx ) \, d\bfx &=& \int _{\tf }   \underset{\bfz \in\tf}{\sup}
 |\STFT_{\Psi}\sigma(\bfz,\mathcal{J}\bfx)| v ( \mathcal{J}\inv
 \mathcal{J}\bfx ) \, d\bfx \\ 
&=& \int _{\hatg \times \cG }   \underset{\bfz \in\tf}{\sup}
 |\STFT_{\Psi}\sigma(\bfz, \bfom )|\, v( \mathcal{J}\inv
  \bfom ) \, d\bfom = \|\sigma \|_{M^{\infty , 1} _{v\circ
     \mathcal{J}\inv} } \, . 
\end{eqnarray*}

\noindent
$(ii) \Rightarrow  (i)$. For the converse,  we note that if $(\bfz ,
\bfom ) = \big( (x,\eta ), \mathcal{J}(\bfy - \bfx ) \big)$ for $\bfz
= (z,\zeta )$ and $\bfom = (\omega , u)$, then $\bfy = (u+z, \zeta )$
and $ \bfx =  (z, \omega + \zeta )$. 
Thus 
$$
| \STFT_{R(g,g)} \sigma (\bfz, \bfom )  |  = |\langle \kn
\pi(z,\omega +\zeta )g, \pi(u+z,\zeta ) g \rangle| \, .
$$
If (ii) holds, then 
$$
\sup _{\bfz \in \tf } |V_\Phi \sigma (\bfz , \bfom ) | \leq
H(\mathcal{J}\inv \bfom )\, ,
$$
and thus
\begin{eqnarray*}
\|\sigma \|_{\mif _{v\circ \mathcal{J}\inv }} &=& \int _{\hattf} \sup
  _{\bfz \in \tf } |V_\Phi \sigma (\bfz , \bfom ) | \, v(
  \mathcal{J}\inv \bfom ) \, d \bfom  \\
&\leq & \int _{\hattf } H(\mathcal{J}\inv \bfom ) v(
\mathcal{J}\inv \bfom ) \, d\bfom = \|H\|_{L^1_v (\tf )} \,.  
\end{eqnarray*}

\medskip
\noindent $(i)\Rightarrow (iii)$:
Let $U$  be  a fundamental
domain of $ \Lambda \subseteq \tf $.  
Set $h(\bfm ) := \sup _{\bfu \in U} \sup _{\bfz \in \tf }
|\STFT_{\Psi}\sigma(\bfz,{\cal J} (\bfm  +\bfu   )|$. Since $\sigma \in
M^{\infty ,1} _{v\circ \cJ ^{-1}} (\tf )$ and $\Phi  \in M^1_{1\otimes
  v\circ \cJ ^{-1}} (\tf )$, Theorem~\ref{wcamalg}
and Corollary~\ref{amalcor} apply and warrant that $h\in \ell ^1_v (\Lambda
)$. 

Next  we use Lemma~\ref{key}(ii) and argue 
as in \eqref{eststft} to obtain 
\begin{equation}
|\langle \kn  \pi(\bfn )g, \pi(\bfm ) g \rangle|    \le 
\underset{\bfz \in\tf}{\sup}|\STFT_{\Psi}\sigma(\bfz,{\cal
  J}(\bfm-\bfn ))| = h(\bfm - \bfn  )
\label{eststft3}
\end{equation}
for $\bfm , \bfn  \in \Lambda $. Thus (iii) is proved. 

\medskip
\noindent $(iii)\Rightarrow (ii)$:
Since  $\piframe$ is a tight frame for $\Ltsp(\tf)$ we 
can express an arbitrary  time-frequency shift $\pi(\bfu)g$ as
\begin{equation}
\pi(\bfu)g = \sum_{\bfm \in \Lambda}\langle \pi(\bfu)g,\pi(\bfm)g \rangle
\pi(\bfm)g.
\label{tfframe}
\end{equation}
 By assumption $g\in\modiiv$ and therefore Theorem~\ref{wcamalg}  implies that $\STFT_g g \in
W(C,\lsp^1_v)(\tf)$. This means that for every relatively compact fundamental
domain $U$ of $\Lambda$ and
\begin{equation}
\alpha(\bfn) = \underset{\bfu \in U}{\sup}|\STFT_g g(\bfn -\bfu)| =
\underset{\bfu \in U}{\sup}|\langle\pi(\bfu)g,\pi(\bfn)g\rangle|,
\quad \bfn \in\Lambda,
\end{equation}
the sequence $\alpha=\{\alpha(\bfn)\}_{\bfn \in \Lambda}$ belongs 
to $\lsp^1_v(\Lambda)$.

Given $\bfy,\bfz \in \tf$ we can write them uniquely as
$\bfy = \bfn +\bfu, \bfz =\bfn'+\bfu'$ for $\bfn,\bfn' \in\Lambda$ and
$\bfu,\bfu' \in U$. Inserting the expansion~\eqref{tfframe} and
the definition of $\alpha$ in the matrix entries we obtain
\begin{align*}
&|\langle K_\sigma \pi(\bfn'+\bfu')g,\pi(\bfn+\bfu)g \rangle| =
|\langle K_\sigma \pi(\bfn')\pi(\bfu')g,\pi(\bfn)\pi(\bfu)g \rangle|\\
&\le \sum_{\bfm,\bfm'\in\Lambda} 
|\langle K_\sigma \pi(\bfn'+\bfm')g,\pi(\bfn+\bfm)g \rangle|
|\langle \pi(\bfu')g,\pi(\bfm')g \rangle|
|\langle \pi(\bfu)g,\pi(\bfm)g \rangle| \\
&\le \sum_{\bfm,\bfm'\in\Lambda} h(\bfn+\bfm-\bfn'-\bfm')
                                 \alpha(\bfm') \alpha(\bfm) \\
& = (h\ast \alpha \ast \tilde{\alpha}) (\bfn-\bfn')
\end{align*}
with $\tilde{\alpha}(\bfn) = \alpha(-\bfn)$. Since $h \in
\lsp^1_v(\Lambda )$
by hypothesis~(iii) and $\alpha \in \lsp^1_v (\Lambda)$ by construction, we
also have $h\ast \alpha \ast \tilde{\alpha} \in \lsp^1_v(\Lambda)$.

Now set
$$H(\bfz) = \sum_{\bfn\in \Lambda} (h\ast \alpha \ast \tilde{\alpha})(\bfn)
\chi_{U-U}(\bfz-\bfv) \qquad \bfz \in \tf \, .$$
Since $\|T_{\bfn } \chi _{U-U}  \|_{L^1_v} \leq v(\bfn ) \|\chi
_{U-U}\|_{L^1_v} $, we obtain that 
$$\|H\|_{\Lsp^1_v}\le \sum_{\bfn\in \Lambda} 
(h\ast \alpha \ast \tilde{\alpha})(\bfn)v(\bfn) \|\chi_{U-U}\|_{\Lsp^1_v}
=c \|h\ast \alpha \ast \tilde{\alpha}\|_{\lsp^1_v} < \infty.$$
For  $\bfy, \bfz \in \tf$ write  $\bfy = \bfn+\bfu$ and $\bfz=\bfn'+\bfu'$ as
before, then we have $\bfy-\bfz \in \bfn-\bfn'+U-U$ and 
$(h\ast \alpha \ast \tilde{\alpha})(\bfn-\bfn') \le H(\bfy-\bfz)$.
Combining these observations we have shown that
$$|\langle K_\sigma \pi(\bfz)g,\pi(\bfy)g \rangle| \le 
(h\ast \alpha \ast \tilde{\alpha})(\bfn-\bfn') \le H(\bfy-\bfz)\, ,$$
and this is (ii).
\end{proof}

\begin{remark}
  We have proved a bit more. The equivalence $ (i) \Leftrightarrow
  (ii) $ requires only that $g\in \modiiv $ without any restriction;
  the implication $(ii) 
  \Rightarrow (iii) $ holds for arbitrary Gabor systems with $g\in
  \modiiv $. Only the implication $(iii) \Rightarrow (ii)$ requires that
  the Gabor system is a frame for $L^2 (\cG )$. 
\end{remark}

Since  the almost diagonalization of implication (i) $\Rightarrow $
(iii)  is important in several applications
(e.g., cf.~\cite{Str05}), we state it explicitly.

\begin{corollary}
  If $\sigma \in  \mod_{v\circ \cJ^{-1}}(\tf)$, $g\in M^1_v (\cG )$ and $\Lambda
  \subseteq \tf $ a lattice in $\tf$, then there   exists a sequence $h
  \in \lsp^1_v( \Lambda)$ such that 
$$
|\langle \kn \pi(\bfm)g ,\pi(\bfn)g \rangle| \le h(\bfm-\bfn),
\qquad \text{for all }\bfm,\bfn \in \Lambda \, . 
$$
\end{corollary}

Next we formulate a similar result on almost diagonalization for
arbitrary LCA groups, even when they do not contain a lattice. Once
more,  we  take recourse to structure theory. 
Recall that  $\tf \simeq \rdd \times \cG _0 \times \widehat{\cG _0}$
with $\cG _0 
\times \widehat{\cG _0}$ containing the compact-open subgroup $\Ksub
\times \Ksub ^\perp $, and let $\bfx \to \dot{\bfx } $ be the canonical
projection from $\tf $ onto $\tf / (\{0\} \times \Ksub \times \Ksub
^\perp ) \simeq \rdd \times \cG _0/\Ksub \times \widehat{\cG _0} / \Ksub
^\perp $. Now let $\Lambda = A\zdd \times D_1 \times D_2 $ be a
quasi-lattice in $\tf $, where $D_1$ is a set of representatives of
$\cG _0/\Ksub $ and $D_2$ a  set of representatives of
$\widehat{\cG _0}/\Ksub ^\perp  $. Then by definition the projection
of $\Lambda$ in $\tf / (\{0\} \times \Ksub \times \Ksub
^\perp )$ is exactly $ \dot{\Lambda } = A\zdd \times \cG _0 / \Ksub
\times \hat{\cG _0}/ \Ksub ^\perp $. Thus $\ind := \dot{\Lambda }$ is a
discrete abelian group. This is the correct index set for the
formulation of the almost diagonalization in general LCA groups. 

Finally, if $v$ is a submultiplicative  weight on $\tf $, then the weight
$\tilde{v}(\dot{\bfx}) =  $ $\sup _{\bfu \in \{0\} \times \Ksub \times
  \Ksub ^\perp } v(\bfx + \bfu )$ is  submultiplicative  on the
quotient $ \tf / (\{0\} \times \Ksub \times \Ksub
^\perp ) $, and $\tilde{v}$ satisfies the GRS-condition if and only if
$v $ does.

\begin{theorem}
\label{th:approxdiag2}
Let $g \in \modiiv(\group)$ for some  admissible weight $v$ on $\tf $, $\Lambda
\subseteq \tf $ be a quasi-lattice, 
and assume that $\piframe$ is a tight Gabor 
frame for $\Ltsp(\group)$. 
Then for $\sigma \in  M^{\infty} (\tf )$ 
the  following properties are equivalent.
\begin{itemize}
\item[(i)] $\sigma \in \mod_{v\circ \cJ^{-1}}(\tf)$.
\item[(ii)] There exists a function $H \in \Lsp^1_v(\tf)$ such that
\begin{equation}
|\langle \kn \pi(\bfz)g,\pi(\bfy)g \rangle| \le H(\bfy - \bfz ),
\qquad \text{for all $\bfy,\bfz \in \tf$.}
\label{}
\end{equation}
\item[(iii)]
There exists a sequence $h_0 \in \lsp^1_v( \ind )$ such that
\begin{equation}
|\langle \kn \pi(\bfn)g ,\pi(\bfm)g \rangle| \le h_0(\dot{\bfm
}-\dot{\bfn } ),
\qquad \text{for all $\bfm,\bfn \in \Lambda$.}
\label{}
\end{equation}
\end{itemize}
\end{theorem}

\begin{proof}
 The equivalence (i) $\Leftrightarrow $ (ii) does not make
reference to any  Gabor frame, and  we have  proved it already  in
  Theorem~\ref{th:approxdiag}. 

\noindent $(i)\Rightarrow (iii)$: 
Let $U= A[0,1)^{2d} \times \Ksub \times  \Ksub ^\perp $  be  a fundamental
domain of $ \Lambda \subseteq \hattf $.  
Set $h(\bfm ) := \sup _{\bfu \in U} \sup _{\bfz \in \hattf }
|\STFT_{\Psi}\sigma(\bfz,{\cal J}(\bfm  +\bfu )  )|$ for $\bfm \in \Lambda$. 
Since $\sigma \in M^{\infty ,1} _{v\circ \cJ ^{-1}} (\tf )$ and $\Phi
\in M^1_{1\otimes v\circ \cJ ^{-1}} (\tf )$,
Corollary~\ref{amalcor} implies  that $h\in \ell ^1_v
(\Lambda )$. Since $\{0\} \times \Ksub \times \Ksub ^\perp $ is a
subgroup of $\tf $, we find that $h(\bfm + \bfu ) = h(\bfm )$ for all
$\bfu \in \{0\} \times \Ksub \times \Ksub ^\perp $. Consequently we
may define  a function $h_0$ on $\ind = \dot{\Lambda }$ unambiguously
by $h_0 (\dot{\bfm }) = h(\bfm )$ for $\dot{\bfm } \in \ind $. Since
 $h\in \ell ^1_v (\Lambda )$, we have $h_0 \in \ell ^1_{\tilde{v}}
 (\ind )$.  

Now we argue as above and   we use Lemma~\ref{key}(ii) and
\eqref{eststft} to obtain  
\begin{equation}
|\langle \kn  \pi(\bfn )g, \pi(\bfm) g \rangle|  \leq  
\underset{\bfz \in\tf}{\sup}|\STFT_{\Psi}\sigma(\bfz,{\cal
  J}(\bfm-\bfn ))| \leq  h(\bfm - \bfn  ) = h_0 (\dot{\bfm} -\dot{\bfn})
\label{eststft2}
\end{equation}
for $\bfm , \bfn \in \Lambda $. Thus (iii) is proved.

\noindent $(iii)\Rightarrow (ii)$:
As in the proof of Theorem~\ref{th:approxdiag} we  express an
arbitrary  time-frequency shift $\pi(\bfu)g$ as 
\begin{equation}
\pi(\bfu)g = \sum_{\bfm \in \Lambda}\langle \pi(\bfu)g,\pi(\bfm)g \rangle
\pi(\bfm)g
\label{tfframe1}
\end{equation}
with respect to a tight Gabor frame $\piframe $. 
The  assumption $g\in\modiiv$ and  Theorem~\ref{wcamalg}
 imply  that $\STFT_g g \in 
W(C,\lsp^1_v)(\tf)$. This means that for the  fundamental
domain $U= A[0,1)^{2d} \times \Ksub \times \Ksub ^\perp $ of $\Lambda$, the
sequence with entries 
\begin{equation}
\alpha(\bfn) = \underset{\bfu \in U}{\sup}|\STFT_g g(\bfn -\bfu)| =
\underset{\bfu \in U}{\sup}|\langle\pi(\bfu)g,\pi(\bfn)g\rangle|,
\quad \bfn \in\Lambda,
\end{equation}
 belongs to $\lsp^1_v(\Lambda)$. As above we note  that $\alpha (\bfn +
 \bfu ) = \alpha (\bfn )$ for $\bfn \in \Lambda $ and $\bfu \in 
\{0\}\times \Ksub \times \Ksub ^\perp $. Thus $\alpha $ can be
identified with a sequence $\alpha _0 (\dot{\bfn}) = \alpha (\bfn
)$ on $\ind $, and $\alpha _0 \in \ell ^1 _{\tilde{v}} (\ind )$. 

Now we follow the proof of Theorem~\ref{th:approxdiag}.
We write  $\bfy,\bfz \in \tf$ in a unique form  as
$\bfy = \bfn +\bfu, \bfz =\bfn'+\bfu'$ for $\bfn,\bfn' \in\Lambda$ and
$\bfu,\bfu' \in U$. Inserting the expansion~\eqref{tfframe1} and
the definition of $\alpha _0$ in the matrix entries we obtain
\begin{align*}
&|\langle K_\sigma \pi(\bfn'+\bfu')g,\pi(\bfn+\bfu)g \rangle| =
|\langle K_\sigma \pi(\bfn')\pi(\bfu')g,\pi(\bfn)\pi(\bfu)g \rangle|\\
&\le \sum_{\bfm,\bfm'\in\Lambda} 
|\langle K_\sigma \pi(\bfn'+\bfm')g,\pi(\bfn+\bfm)g \rangle|
|\langle \pi(\bfu')g,\pi(\bfm')g \rangle|
|\langle \pi(\bfu)g,\pi(\bfm)g \rangle| \\
&\le \sum_{\bfm,\bfm'\in\Lambda} h_0(\dot{\bfn }+\dot{\bfm
}-\dot{\bfn' }-\dot{\bfm' }) \alpha _0(\dot{\bfm' }) \alpha
_0(\dot{\bfm}) \\ 
 = & (h_0\ast \alpha_0 \ast \widetilde{\alpha _0}) (\dot{\bfn }-\dot{\bfn'
   })
\end{align*}
with $\widetilde{\alpha _0}(\bfn) = \alpha _0(-\dot{\bfn })$. Here it is
crucial that $\ind = \dot{\Lambda }$ is a group.  Since $h_0 \in
\lsp^1_{\tilde{v}} (\ind ) $
by hypothesis~(iii) and $\alpha _0 \in \lsp^1_{\tilde{v}}(\ind )$ by
construction, the sequence  $h_0\ast \alpha _0 \ast \tilde{\alpha _0} $ is
also in $\lsp^1_{\tilde{v}}(\ind)$.

Now set
$$H(\bfz) = \sum_{\bfn\in \Lambda} (h_0\ast \alpha _0 \ast
\tilde{\alpha _0})(\dot{\bfn })
\chi_{U-U}(\bfz-\bfn).$$
Since $\|T_{\bfn } \chi _{U-U}  \|_{L^1_v} \leq v(\bfn ) \|\chi
_{U-U}\|_{L^1_v} $ and $\tilde{v}(\dot{\bfn }) \leq C v(\bfn )$, we obtain that 
$$\|H\|_{\Lsp^1_v}\le \sum_{\bfn\in \Lambda} 
(h_0 \ast \alpha _0 \ast \tilde{\alpha _0})(\dot{\bfn })v(\bfn)
\|\chi_{U-U}\|_{\Lsp^1_v} 
=c \|h_0\ast \alpha_0 \ast \tilde{\alpha_0}\|_{\lsp^1_{\tilde{v}}} < \infty.$$
Arguing as before, we  show that 
$$|\langle K_\sigma \pi(\bfz)g,\pi(\bfy)g \rangle| \le 
(h_0\ast \alpha _0\ast \widetilde{\alpha _0})(\dot{\bfn }-\dot{ \bfn'} )
\le H(\bfy-\bfz)\, ,$$ 
and this is (ii).
\end{proof}

Remark: Despite the resemblance of the proofs,
Theorem~\ref{th:approxdiag} is not a special case of 
Theorem~\ref{th:approxdiag2} because in the former case we consider
an arbitrary lattice in $\tf$, if it exists, whereas in the latter case
we consider a very special quasi-lattice that respects the
factorization of $\tf $ as $\rdd \times \cG _0 \times \hat{\cG _0}$.

\section{Sj\"ostrand's Results on Locally Compact Abelian Groups}
\label{s:sjolca}

The characterization of almost diagonalization through time-frequency
properties of the symbol leads to the generalization of Sj\"ostrand's
results to LCA groups. Note that in no place do we resort to typical
arguments from pseudodifferential operator calculus. 

\medskip
\noindent
{\bf Boundedness on $\Ltsp(\group)$:} First we prove that any
pseudodifferential operator $K_\sigma $ with a symbol in the
generalized Sj\"ostrand class is bounded  on all modulation spaces
with appropriate weight. As a preparation we need a lemma on the
properties of the Rihaczek distribution $ R(f,g)(x,\xi )  = f(x)
\overline{\hat{g}(\xi )} \,\, \overline{\langle \xi, x\rangle}$
that generalizes~\cite{Gro04b} to     LCA groups and weighted
modulation spaces.

\begin{lemma}
\label{wdstft}
(i)   Let  $\varphi , \psi , f,g \in L^2(\cG ) $ and set $\Phi=R(\varphi , \psi
  ) \in L^2 (\tf ) $. Then, with $\bfx =(x,\xi)\in \tf $,
  $\bfom =(\omega , u)\in \hattf $, we have
\begin{equation}\label{1433}
\cV_\Phi\bigl(R(g,f)\bigr) (\bfx, \bfom )=\overline{ \langle \xi , u
  \rangle } \,
V_\psi  g(x, \xi+\omega ) \, \overline{V_\psi f(x+ u,  \xi  )} \, .
\end{equation}
(ii) If $f\in M^{p,q}_{m} (\cG )$ and $g \in M^{p', q' } _{1/m} (\cG )$,
 then $R(g,f) \in M^{1,\infty } _{1/v\circ \cJ \inv }
(\tf )$ and 
\begin{equation}
  \label{eq:chr4}
  \| R(g,f) \|_{ M^{1,\infty } _{1/v\circ \cJ \inv }} \leq C
    \|f\|_{M^{p,q}_{m}} \, \|g\|_{M^{p', q' } _{1/m} } \, .
\end{equation}
\end{lemma}

\begin{proof}
 The proof is similar to~\cite{Gro04b}. We write  the time-frequency
 shifts of the \riha\ distribution  explicitly as 
  \begin{equation*}
        M_\bfom  T_\bfx  R(\varphi ,\psi )(t,\tau ) = \langle \omega
        ,t\rangle \langle \tau , u\rangle  \, \varphi (t-x) 
\, \overline{\hat{\psi} (\tau -\xi) \, \langle \tau -\xi, t-x\rangle }
\, .
  \end{equation*}
Consequently, after a substitution, 
\begin{eqnarray*}
  \cV _\Phi R(g,f) (\bfx ,\bfom ) &=& \langle R(g,f), M_\bfom  T_\bfx R(\varphi ,
  \psi )\rangle \\
&=& \iint _{\tf } g(t) \overline{\hat{f}(\tau )} \, \overline{\langle \tau
,t\rangle }  \, \overline{\varphi (t-x)} \, \hat{\psi} (\tau -\xi)\,  
\langle \tau -\xi,   t-x \rangle  \, \overline{\langle \omega
  ,t \rangle \langle \tau ,u\rangle} \,\, dtd\tau \\
&=&  \langle \xi , x\rangle \, \int _{\cG } g(t)  \overline{\varphi
  (t-x)} \, \overline{\langle \xi +\omega , t\rangle}\, dt \,\, \cdot \,
 \int _{\hatg } \overline{\hat{f}(\tau )} \,\hat{\psi} (\tau -\xi)\, 
\langle \tau , -x-u\rangle  \, d\tau  \\
&=&  \langle \xi ,x \rangle \, V_{\varphi } g (x, \xi+\omega) \,
\overline{V_{\hat{\psi } }\hat{f} (\xi, -x-u)} \\
&=& \overline{\langle \xi  , u\rangle } \, V_{\varphi } g (x, \xi+\omega) \,
\overline{V_{\psi } f (x+u, \xi  )} \, . 
\end{eqnarray*}
In the last transformation  we have used the fundamental
formula~\eqref{STFTft}. Since both $R(g,f)$ and $R(\varphi, \psi )$ are in
$L^2(\tf  )$, the integral defining $\cV _\Phi R(g,f)$ is absolutely
convergent on $\tf $, and so the application of Fubini's theorem is
justified. 

(b) is a consequence of (a). For simplicity we use the window $\Phi =
R(\varphi , \varphi )$ and use the fact that different windows in
$M^1_v(\cG )$ yield equivalent norms on $M^{p,q}_m (\cG )$
(Lemma~\ref{equivnorm}). Consequently,   
\begin{eqnarray*}
  \| R(g,f)\|_{ M^{1,\infty } _{1/v\circ \cJ \inv }} &=& \sup _{\bfom
    \in \hattf } \frac{1}{v(\cJ \inv \bfom )} \, \int _{\tf } |V_\Phi
  (R(g,f)(\bfx , \bfom )| \, d\bfx \\
&=& \sup _{(\omega ,u) \in \hattf } \frac{1}{v(u,-\omega )} \, \iint _{\tf} |V_{\varphi
  } f (x+u, \xi ) | \, |V_\varphi g ( x, \xi +\omega )| \, dx d\xi =
  (\ast ) \
\end{eqnarray*}
Since $m(x,\xi +\omega ) \leq C v(-u,\omega ) m(x+u, \xi)$ 
by~\eqref{moderateweights} and since $v$ is even, we can continue 
the estimate by 
\begin{eqnarray*}
  (\ast ) &\leq & C \sup _{(\omega ,u) \in \hattf } \iint _{\tf }  |V_{\varphi
  } f (x+u, \xi ) | m(x+u,\xi )  \, |V_\varphi g ( x, \xi +\omega )|
  \frac{1}{m(x, \xi +\omega )} \, dx d\xi  \\
&\leq & C \sup _{(\omega ,u) \in \hattf } \| V_{\varphi } f\,  m
\|_{L^{p,q}} \, \| V_{\varphi } g \, m^{-1} \|_{L^{p',q'}} \\
&=& C \|f\|_{M^{p,q}_{m}} \, \| g\|_{M^{p',q'}_{1/m}} \, ,
\end{eqnarray*}
where in the last step we have applied H\"older's inequality. 
 \end{proof}

We are now ready to prove that operators in the Sj\"ostrand class are
bounded on modulation spaces. 

\begin{theorem}
  \label{boundedness}
Let   $v$ be  an admissible weight on $\tf $. If $\sigma\in \modvj(\tf)$, then
$K_\sigma $ is bounded on all modulation spaces $M^{p,q}_m(\cG )$ for
$1\leq p,q \leq \infty $  and every $v$-moderate weight $m$. In
particular, $\kn $ is bounded on $L^2(\cG )$. 
\end{theorem}

\begin{proof}
We apply the duality~\eqref{eq:chr5} and Lemma~\ref{wdstft}(ii) to   obtain 
\begin{eqnarray*}
|\langle \kn  f,g\rangle | &=& |\langle \sigma , R(g,f)\rangle | \\
&\leq &  C\, \|\sigma \|_{\mif _{v\circ \cJi }} \,  \| R(g,f) \|_{M^{1,\infty
    }_{1/v\circ \cJi}} \\
&\leq &   C' \, \|\sigma \|_{\mif _{v\circ \cJi }} \,   \|f\|_{M^{p,q}_{m}} \,
  \|g\|_{M^{p', q' } _{1/m} } \, . 
  \end{eqnarray*}
Since this inequality holds for \emph{all} $g \in M^{p',q'}_{1/m}(\cG
)= (M^{p,q}_m(\cG))^*$, we have shown that 
$$
\|\kn  f \|_{M^{p,q}_m} \leq C\,  \|\sigma \|_{M^{\infty ,1}_{v\circ \cJi}} \,
\|f\|_{M^{p,q}_m} \, .
$$
If $(p,q) = 
(1,\infty )$ or $(p,q) = (\infty ,1 )$, we observe that these spaces
are also dual spaces of a modulation space~\cite{BGH05}, thus we have
proved the boundedness of $\kn $ for all parameters $p,q \in [1,\infty
]$. For $p=q=2$ and $m\equiv 1$, we obtain the boundedness on
$M^{2,2}(\group )
= L^2(\cG )$. 
 \end{proof}

\medskip
\noindent
{\bf The Banach Algebra Property.} 
Whereas the boundedness property uses typical arguments from
time-frequency analysis, the Banach algebra property lies much deeper
and requires the characterization of the generalized Sj\"ostrand class
via almost diagonalization. 

Recall that $\cC _v (\ind) $ is
the Banach algebra of all matrices on the index set $\ind $ that
are dominated by convolution operators in $\ell ^1_v (\ind )$. 

\begin{theorem}
\label{th:banach}
The space $\cA = \Op(\modvj(\tf))$ is a Banach algebra with respect to the composition
of operators and with the norm $\|K_\sigma \|_{\cA } := \|M(\sigma
)\|_{\cC _v} \asymp 
\|\sigma \|_{\mif _{v\circ \mathcal{J}\inv}}$.
\end{theorem}
\begin{proof}
Let $\gabframe$ be a tight Gabor frame with $g\in \modiiv (\cG )$ with
respect to a quasi-lattice $\Lambda \subseteq \tf $. Recall that the
projection of $\Lambda $ into a quotient of $\tf $ results in the
discrete abelian group $\ind $.  We may  
assume without loss of generality that $g$ is normalized such that its
(lower and upper) frame bound is 1.
Using~Lemma~\ref{twistedlemma} we compute
\begin{gather}
M({\cal F}^{-1}(\hat{\sigma} \twist \hat{\tau})) C_g f =
C_g (\knst f) = C_g (\kn \knt f) \notag \\
 = M(\sigma) C_g (\knt f) = M(\sigma) M(\tau) C_g f.
\label{}
\end{gather}
This means that  $M({\cal F}^{-1}(\hat{\sigma} \twist \hat{\tau}))$ and $M(\sigma) M(\tau)$
coincide  on $\ran(C_g)$. 
Since $\kn, \knt$, and  $ \knst$ are all bounded on $\Ltsp(\group)$,
it follows  
from Lemma~\ref{le:range} that $\ker(M({\cal F}^{-1}(\hat{\sigma} 
\twist \hat{\tau})))= \ker(M(\sigma) M(\tau))$. Hence
\begin{equation}
M({\cal F}^{-1}(\hat{\sigma} \twist \hat{\tau}))= M(\sigma) M(\tau)
\label{}
\end{equation}
on $\ltsp(\Lambda)$.

Now let $\sigma , \tau \in \mif _{v \circ \cJi} (\tf )$. Then, by
Theorem~\ref{th:approxdiag2} 
$M(\sigma ) \in \cC _v (\ind )$ and $M(\tau ) \in \cC _v (\ind)
$. Since $\cC _v (\ind ) $ is a Banach algebra, the product $M(\sigma
) M(\tau )$ is also in $\cC _v (\ind )$, and furthermore  
$$
\| K_{{\cal F}^{-1}(\hat{\sigma} \twist \hat{\tau})} \|_{\cA } 
= \|M({\cal F}^{-1}(\hat{\sigma} \twist \hat{\tau}))\|_{\cC _v} \leq
\|M(\sigma )\|_{\cC _v} \|M(\tau )\|_{\cC _v} =  
\|K_\sigma \|_{\cA } \|K_\tau \|_{\cA } \, .
$$
 By Theorem~\ref{th:approxdiag2} we conclude that 
${\cal F}^{-1}(\hat{\sigma} \twist \hat{\tau})) \in \mif _{v \circ
  \cJi }$ and that $K_\sigma K_\tau \in \cA  = \Op(\modvj(\tf))$. 
\end{proof}


\medskip
\noindent
{\bf The Wiener Property:}
We now state the main result of this paper, the Wiener property of the
generalized Sj\"ostrand class. This is the deepest theorem of this
paper and requires the combination of all methods developed so far,
namely the almost diagonalization, the Wiener property of the matrix
algebra $\cC _v $, and the existence and  properties of tight  Gabor frames. 

\begin{theorem}
\label{mainwiener}
Let $v$ be an admissible weight. If $\sigma\in \modvj(\tf)$ and
if $\kn$ is invertible on $\Ltsp(\group)$, then $(\kn)^{-1}  = K_\tau$
for some $\tau\in \modvj(\tf)$.
\end{theorem}
\begin{proof}
As in the proof of Theorem~\ref{th:banach} we use a tight Gabor frame 
$\gabframe$ with $g\in \modiiv$ and 
with (lower and upper) frame bounds equal to 1. Let $\tau$ be the unique distribution
such that $(\kn)^{-1} = K_\tau$. By Lemma~\ref{le:range} we have that
the matrix $M(\tau,g,\Lambda)=M(\tau)$ is bounded on $\ltsp(\Lambda)$ and 
maps $\ran(C_g)$ into $\ran(C_g)$ with 
$\ker(C^{\ast}_g) \subseteq \ker(M(\tau))$. 

We show that $M(\tau)$ is the pseudoinverse of $M(\sigma)$. Let 
$c = C_g f \in \ran(C_g)$, then
\begin{equation}
\label{}
M(\tau) M(\sigma) C_g f = M(\tau) C_g (\kn f ) = C_g (K_\tau \kn f)
= C_g f,\label{matcalc}
\end{equation}
where we have used~\eqref{matrix2} and the property that $\gabframe$
is a tight frame. Relation~\eqref{matcalc} says that
$M(\tau) M(\sigma) = I$ on $\ran(C_g)$. Furthermore,
$\ker(M(\sigma)), $ $\ker(M(\tau)) \supseteq \ran(C_g)^\perp$, thus we conclude
that $M(\tau)$ is the pseudoinverse of $M(\sigma)$.

By Theorem~\ref{th:approxdiag2} the property $\sigma\in \modvj(\tf)$ implies 
that $M(\sigma) \in \nswv(\ind)$. Applying
Corollary~\ref{cor:pinv} to $M(\sigma )$,  we deduce that  
$M(\tau) = M(\sigma)^{+} \in \nswv(\ind)$. Using
Theorem~\ref{th:approxdiag2} once more, we have  shown  that $K_\tau
\in \modvj(\tf)$. 
\end{proof}

\medskip
\noindent
{\bf  Spectral invariance on modulation spaces:} According to 
Theorem~\ref{mainwiener} the inverse $K_\sigma \inv $ has again a
symbol in $\modvj(\tf)$. Consequently, by Theorem~\ref{boundedness}
$K_\sigma \inv $ acts boundedly on a large class of modulation spaces
depending only on the class of the weight $v$. 

\begin{corollary}
  Let $v$ be an admissible weight. If $\sigma\in \modvj(\tf)$ and
if $\kn$ is invertible on $\Ltsp(\group)$, then $\kn $ is invertible
simultaneously on all modulation spaces $M^{p,q}_m(\cG )$ for $1\leq
p,q \leq \infty $ and $v$-moderate weight $m$. 
\end{corollary}

\section{Special Groups} \label{s:special}

\subsection{Sj\"ostrand's Class} \label{ss:sjoclass}

Clearly, when we choose $\group=\Rst ^d$ (whence $\hatg =
\widehat{\Rst ^d}=\Rst ^d$)
in the derivations of the previous sections, then we recover
Sj\"ostrand's results~\cite{Sjo94,Sjo95}. The additional insight
gained from our approach is the identification of the cornerstones of
Sj\"ostrand's result, namely modulation spaces and the corresponding
time-frequency techniques, the striking appearance of certain matrix
algebras and their spectral invariance, and the almost diagonalization
by Gabor frames. 

\subsection{Discrete Pseudodifferential Operators} \label{ss:discrete}

Let us consider the case $\group = \Zst$, $\hatg = \Tst$. Thus $\kn$,
now acting on $\ltsp(\Zst)$, becomes a discrete pseudodifferential
operator. We refer to~\cite{Mas87} for a detailed review of discrete 
pseudodifferential operators.
Of course, the action of $K_\sigma $ can be described simply by a
matrix. The next lemma elucidates the relation between the symbol
class $\modvj(\Zst\times\Tst)$, the ``discrete Sj\"ostrand class'',
and the corresponding class of matrices. A similar calculation was
made in~\cite{Gal00}. 


\begin{lemma}
Let $\kn$ be a pseudodifferential operator defined on $\ltsp(\Zst)$ and
let $v$ be a weight function on $\Zst$. Then the matrix corresponding
to  $\kn $ is in $ \nswv(\Zst)$ 
if and only if $\tii\sigma \in \modvij(\Zst\times\Tst)$, where
$\tii\sigma(x,\xi):=\sigma(x,-\xi)$.
\end{lemma}
In this special case Theorem~\ref{mainwiener} coincides with Theorem~\ref{th:bask} of
Gohberg et al.\ and Baskakov. However, our proof of  Theorem~\ref{mainwiener} does not give a
new proof of Theorem~\ref{th:bask}, because we have used the Gohberg-Baskakov
result in the proof of the general theorem. 

\vspace{3 mm}

\begin{proof}
We begin by calculating $\|\tii\sigma\|_{\modvj(\Zst\times\Tst)}$.
We first compute $\STFT_{\Psi}\sigma$ for an appropriate window $\Psi$
in $\modiivi(\Zst \times \Tst)$.
We choose $\Psi = \delta \otimes {\bf 1}$. Since $\Psi $ is the
characteristic function of the compact-open subgroup $\{0\} \times
\Tst$ in $\Zst\times\Tst$, its STFT $V_\Psi \Psi $ is integrable with
respect to every submultiplicative weight $v$ by~\eqref{eq:chr1}, and thus $\Psi
\in M^{1}_{v\otimes 1}(\Zst\times\Tst)$ for every submultiplicative
$v$ on $\Zst $. 

Let $\bfx =(x,\xi) \in \Zst\times\Tst , \bfom = (\omega,u)\in
\Tst \times   \Zst$. We compute
\begin{align*}
\STFT_{\Psi}\tii\sigma(\bfx,\bfom) & = \sum_{z\in \Zst}\int \limits_{\Tst}  
\tii \sigma(z,\zeta) \overline{ M_{\bfom} T_{\bfx} \Psi(z,\zeta)} \,d\zeta \\ 
& = \sum_{z\in \Zst}\int \limits_{\Tst}  
\sigma(z,-\zeta) \delta(z-x) {\bf 1}(\zeta-\xi)
e^{-2\pi i \zeta u} e^{-2\pi i z \omega} \,d\zeta \\ 
& =  e^{-2\pi i x\omega } \, \int \limits_{\Tst} \sigma(x,-\zeta)
e^{-2\pi i \zeta u} \,d\zeta \\ 
& = e^{-2\pi i x \omega } {\cal F}_2^{-1} \sigma(x,u).
\end{align*}
Consequently, since $(v\otimes 1)( \cJi (\omega ,u)) = (v\otimes
1) (u, -\omega ) = v(u)$, we obtain 
\begin{align}
\|\tii\sigma\|_{\modvij} & = \sum_{u\in \Zst}\int \limits_{\Tst} 
\underset{\bfx}{\sup} |\STFT_{\Psi}\tii \sigma(\bfx,\bfom)| \,  
(v\otimes 1)( \cJi)(\omega,u)) \,d\omega
\notag \\
& = \sum_{u\in \Zst}\int \limits_{\Tst} \underset{x,\xi}{\sup} 
|e^{- 2\pi i x\omega } {\cal F}_2^{-1} \sigma(x,u)|\, v(u) \,d\omega 
\notag \\
& = \sum_{u\in \Zst} \underset{x\in \Zst }{\sup} |{\cal F}_2^{-1} \sigma(x,u)|
\, v(u).
\label{z2}
\end{align}
Next we compute $\kn f(x)$.
\begin{align}
\kn f(x) & = \int \limits_{\Tst} \sigma(x,\xi)\hat{f}(\xi)
             e^{2\pi i x\xi}\,d\xi \notag \\
& = \sum_{u\in\Zst}\Big(\int \limits_{\Tst} \sigma(x,\xi)
e^{-2\pi i (u-x)\xi} d\xi \Big) f(u)  \notag \\
& = \sum_{u\in\Zst} {\cal F}_2 \sigma(x,u-x) f(u) \notag \\
& = \sum_{u\in\Zst} {\cal F}_2^{-1} \sigma(x,x-u) f(u) \notag \\
& = Af \, . \label{z3}
\end{align}
This means that the matrix $A$ corresponding to $K_\sigma $ has the 
 entries 
\begin{equation}
A_{x,u}:={\cal F}_2^{-1}\sigma(x,x-u), \qquad x,u\in \Zst.
\label{Amatrix}
\end{equation}
Comparing~\eqref{z2},~\eqref{z3} and~\eqref{Amatrix}, we find that  the assumption
$\tii \sigma \in \modvij(\Zst\times\Tst)$ yields
$$
\sum_{u\in \Zst} \underset{x\in \Zst}{\sup} |{\cal F}_2\tii \sigma(x,u)| v(u)
< \infty \Longleftrightarrow 
\sum_{u\in \Zst} \underset{x\in \Zst }{\sup} |A_{x,x-u}| v(u) < \infty \, . $$
Comparing with  Definition~\ref{def:nsw}, we have shown that the
matrix $A$ corresponding to $K_\sigma $ is in $\cC _v( \Zst )$. 
\end{proof}

\subsection{Periodic Pseudodifferential Operators} \label{ss:periodic}

We consider periodic pseudodifferential operators,
cf.~e.g.~\cite{hormander3,SV02}. 
In this case $\group = \Tst$, $\hatg = \Zst$ and thus the symbol is 
defined on $\Tst \times \Zst$. Analogous to the previous section we will
first analyze $\modivj(\Tst\times\Zst)$.

We first compute $\STFT_{\Psi}\sigma$ for the  window $\Psi = {\bf 1}
\otimes \delta$. As before  $\Psi
\in \modiiiv(\Tst \times \Zst)$ for any submultiplicative weight on
$\Zst $. 

Let $\bfx = (x,\xi)\in \Tst \times \Zst , \bfom =(\omega,u) \in  \Zst
\times \Tst $. We compute 
\begin{align}
\STFT_{\Psi}\sigma(\bfx,\bfom) & = \int \limits_{\Tst}  \sum_{\zeta\in \Zst}
\sigma(z,\zeta) \overline{ M_{\bfom} T_{\bfx} \Psi(z,\zeta)} \,dz \\ 
& = \int \limits_{\Tst}  \sum_{\zeta\in \Zst}
\sigma(z,\zeta) {\bf 1}(z-x)   \,  \delta(\zeta-\xi) \, 
e^{-2\pi i z \omega} e^{-2\pi i \zeta u} \,dz \\ 
& = e^{-2\pi i \xi u} \, \int \limits_{\Tst} \sigma(z,\xi) e^{-2\pi i z\omega}\,dz \\ 
& = e^{-2\pi i \xi u}  {\cal F}_1 \sigma(\omega,\xi).
\end{align}
Consequently, since $(1\otimes v)(\cJi (\omega, u)) = (1\otimes
v)(u,-\omega ) = v(-\omega )= v(\omega )$, we have
\begin{align}
\|\sigma\|_{\modivj} & = \int \limits_{\Tst} \sum_{\omega\in \Zst}
\underset{\bfx}{\sup} |\STFT_{\Psi} \sigma(\bfx,\bfom)|(1\otimes v)(\cJi(\bfom))
\,du \\
& = \sum_{\omega\in \Zst} \underset{\xi \in \Zst }{\sup} |{\cal F}_1 \sigma(\omega,\xi)|
v(\omega).
\label{z4}
\end{align}
Next we compute the Fourier transform of $\kn f$.
\begin{align}
({\cal F}\kn f)(\xi) & = {\cal F}\Big(\sum_{\omega\in\Zst}\sigma(x,\omega)
\hat{f}(\omega) e^{2\pi i x \omega}\Big)(\xi) \notag \\
& = \int \limits_{\Tst} \sum_{\omega\in\Zst}\sigma(x,\omega)
\hat{f}(\omega) e^{2\pi i x (\omega-\xi)}dx \notag \\
& = \sum_{\omega\in\Zst} {\cal F}_1
\sigma(\xi-\omega,\omega)\hat{f}(\omega)\, .\label{z5}
\end{align}
Let  $A$ be   the matrix  with entries
$A_{\xi , \omega }= {\cal F}_1 \sigma(\xi-\omega,\omega ), \, \xi , \omega
\in \Zst $. Then 
$$
A \hat{f} = \cF K_\sigma \cF \inv \hat{f} \, ,
$$ 
in other words, $A$ describes the action of $K_\sigma $ on the Fourier 
coefficients  $\hat{f} $ of $f$.   Using~\eqref{z4} and~\eqref{z5} we   
   see that 
   \begin{eqnarray*}
    \|A \|_{\cC _v} &=& \sum _{\omega \in \Zst } \sup _{\xi\in \Zst }
    |A _{\xi, \xi -\omega } | v(\omega ) \\
&=& \sum _{\omega \in \Zst } \sup _{ \xi \in \Zst }  | {\cal F}_1
\sigma(\omega,\xi - \omega)| \, v(\omega ) \\
&=& \sum _{\omega \in \Zst } \sup _{ \xi \in \Zst }  | {\cal F}_1
\sigma(\omega,\xi )| \, v(\omega ) \\
&=& \|\sigma \|_{\modivj     } \, . 
   \end{eqnarray*}
So we have shown that 
${\cal F}\kn {\cal F}^{-1} \in \nswv(\Zst)$ if and
only if $\sigma \in \modivj(\Tst\times \Zst)$.

\bibliographystyle{abbrv}

\def\cprime{$'$} \def\cprime{$'$}


\end{document}